\documentclass[12pt]{article}

\usepackage{graphicx}
\usepackage{color}
\usepackage{amsfonts}
\usepackage{amsmath}
\usepackage{multirow}
\usepackage{amssymb}
\usepackage{changes}
\usepackage{apptools}

\begin{document}
 
\newcommand{\al}{\mbox{$\alpha$}}
\newcommand{\be}{\mbox{$\beta$}}
\newcommand{\ep}{\mbox{$\epsilon$}}
\newcommand{\gam}{\mbox{$\gamma$}}
\newcommand{\sig}{\mbox{$\sigma$}}
\newcommand{\sign}{\mbox{sgn}}

\DeclareRobustCommand{\FIN}{%
  \ifmmode 
  \else \leavevmode\unskip\penalty9999 \hbox{}\nobreak\hfill
  \fi
  $\bullet$ \vspace{5mm}}

\newcommand{\calA}{\mbox{${\cal A}$}}
\newcommand{\calB}{\mbox{${\cal B}$}}
\newcommand{\calC}{\mbox{${\cal C}$}}

\newcommand{\muas}{\mbox{$\mu$-a.s.}}
\newcommand{\Nat}{\mbox{$\mathbb{N}$}}
\newcommand{\Rea}{\mbox{$\mathbb{R}$}}
\newcommand{\Prob}{\mbox{$\mathbb{P}$}}

\newcommand{\nin}{\mbox{$n \in \mathbb{N}$}}
\newcommand{\suc}{\mbox{$\{X_{n}\}$}}
\newcommand{\sucP}{\mbox{$\mathbb{P}_{n}\}$}}

\newcommand{\conv}{\rightarrow}
\newcommand{\convn}{\rightarrow_{n\rightarrow \infty}}
\newcommand{\convp}{\stackrel{p}{\rightarrow}}
\newcommand{\convs}{\stackrel{a.s.}{\rightarrow}}
\newcommand{\convw}{\stackrel{w}{\rightarrow}}
\newcommand{\convd}{\stackrel{\cal D}{\rightarrow}}

\newtheorem {Prop}{Proposition} [section]
\newtheorem {Lemm}[Prop] {Lemma}
\newtheorem {Theo}[Prop]{Theorem}
\newtheorem {Coro}[Prop] {Corollary}
\newtheorem {Nota}[Prop]{Remark}
\newtheorem {Ejem}[Prop] {Example}
\newtheorem {Defi}[Prop]{Definition}
\newtheorem {Figu}[Prop]{Figure}
\AtAppendix{\counterwithin{Prop}{subsection}}

\title{\sc The complex behaviour of Galton rank order statistic.\footnote{Research partially supported by 
FEDER, Spanish Ministerio de Econom\'{\i}a y Competitividad, grants MTM2014-56235-C2-1-P and MTM2017-86061-C2-1-P and Junta de Castilla y Le\'on, grants VA005P17 and VA002G18.}}

\author{E. del Barrio$^{1}$, J.A. Cuesta-Albertos$^{2}$\\ and C. Matr\'an$^{1}$ \\
$^{1}$\textit{Departamento de Estad\'{\i}stica e Investigaci\'on Operativa and IMUVA,}\\
\textit{Universidad de Valladolid} \\ $^{2}$ \textit{Departamento de
Matem\'{a}ticas, Estad\'{\i}stica y Computaci\'{o}n,}\\
\textit{Universidad de Cantabria}}
\maketitle

\vspace{-1cm}
\begin{abstract}{
Galton's rank order statistic is one of the oldest statistical tools for two-sample comparisons. It is also a very natural index to measure departures from stochastic dominance. Yet, its asymptotic behaviour has been investigated only partially, under restrictive assumptions. This work provides a comprehensive {study} of this behaviour, based on the analysis of the so-called contact set (a modification of the set in which the quantile functions coincide).
We show that a.s. convergence to the population counterpart holds if and only if {the} contact set has zero Lebesgue measure. When this  set is
finite we show that the asymptotic behaviour is determined by the  local behaviour of a suitable reparameterization  of the quantile functions  in a neighbourhood of the contact points. Regular crossings result in standard rates and Gaussian limiting distributions, but higher order contacts (in the sense introduced in this work) or contacts at the extremes of the supports may result in different rates and non-Gaussian limits.} \end{abstract}

\noindent
{\it Keywords:}  Relaxed stochastic dominance, asymptotics, consistency, Galton rank order statistic, comparison of quantile functions, contact points, crossings, tangencies, contact intensity.


\section{Introduction and main results}

 The Introductory Remarks in Darwin's report on the benefits of cross-fertilization to the 
propagation of vegetal species \cite{Darwin} include the following comment, 
by Galton:
{\it ``The observations\dots have no prim\^{a} facie appearance of regularity. But as 
soon as we arrange them  in order of their magnitudes,\dots {. We now see,} with few exceptions, 
{that}\ldots the largest plant on the crossed side\dots exceeds the largest plant on the self-fertilised 
side,  {that}\ldots the second exceeds the second,\dots and so on\dots "}. {With this argument,} 
Galton opened a simple   way of comparison of  distributions, just by comparing the values with the 
same ranks  in their respective settings. 
 
Given two samples of equal size,  $X_1,\dots,X_n$ and $Y_1,\dots,Y_n$, respectively coming from 
the distribution functions (d.f.'s in the sequel) $F$ and $G$,  let us denote\footnote{We have tried to use throughout standard or natural notation. However, a complete enough notation guide is included   at the end of this section.} by $F_n$ and $G_n$ 
the corresponding sample d.f.'s. Galton's solution consisted in reordering both data samples in 
increasing order:  $X_{(1)},\dots,X_{(n)}$ (coming from the control) and $Y_{(1)},\dots,Y_{(n)}$  
(from the treatment) and computing $\mathcal G(F_n,G_n):=\#\{i: \;  X_{(i)} > Y_{(i)}\}$, concluding 
improvement under the treatment whenever $\mathcal G(F_n,G_n)$  is small enough. When $F=G$ { is continuous,} the 
distribution of this `Galton Rank Order' statistic is uniform on $\{0,1,\dots,n\}$ (see \cite{Chung-Feller};
see also \cite{Sparre-Andersen}, \cite{Hodges-Galton} or \cite{Feller1} for alternative proofs). As explained in   \cite{Hodges-Galton},  in Darwin's problem the 
sample sizes were 15 and $\mathcal G(F_{15},G_{15})=2$, thus the $p$-value associated to Galton's 
approach is 3/16, which is not as rare as he suspected. 

Galton's strategy  was related to the assessment of stochastic dominance of $G$ over $F$, 
$F<_{st}G$, being the alternative to the null hypothesis $F=G$. Recall that,  by definition,
\begin{equation*}
F \leq_{st}G \mbox{ whenever } F(x)\geq G(x) \mbox{  for every } x\in\Rea.
\end{equation*}

As noted in \cite{Lehmann}, this relation is better understood when it is 
stated in terms of the quantile functions: if  
$F^{-1}$ is the  quantile function associated to $F,$ defined by 
\begin{equation}\label{defq}
F^{-1}(t):= \inf\{x : \ t\leq F(x)\}, 
\mbox{ for } t \in (0,1), 
\end{equation}
then
\begin{equation*}
F \leq_{st}G \mbox{   whenever } F^{-1}(t)\leq G^{-1}(t) \mbox{    for every } t\in (0,1).
\end{equation*}

A useful feature of {the} quantile functions is that they provide a canonical representation of {random variables (}r.v.'s {in that follows)} with a 
given d.f.: if we consider {the} Lebesgue measure, $\ell$, on the unit interval $(0,1)$, the  function $F^{-1}$ 
is a r.v. with d.f. $F$.  With this in mind, 
we set 
\begin{equation}\label{gamma}
\gamma(F,G):=\ell\{t  : F^{-1}(t)>G^{-1}(t)\}
\end{equation}
and observe that
\begin{equation}\label{GaltonROS}
\mathcal G(F_n,G_n) =n\gamma(F_n,G_n).
\end{equation}

Early work on Galton's rank statistic focused on the case $F=G$ and equal sample sizes. 
Special mention  should be given to \cite{Csaki}, which analyzes the joint 
behaviour of the Kolmogorov-Smirnov and Galton statistics (under $F=G$).  Also for equal sample sizes, later, \cite{Gross} considered
the intermediate case with $F\ne G$ possibly, but $\ell\left\{F^{-1}=G^{-1}\right\}>0$. 
Focusing on the dominance model $F=G$ vs $F\leq_{st}G$, \cite{BehnenNeuhau1983} addressed the local 
asymptotic efficiency of $\gamma(F_{n},G_m)$, noting that it is just 
a generalization of Galton's statistic (recall \eqref{GaltonROS}) and using
empirical processes techniques  to obtain the asymptotic {distribution} of $\gamma(F_n,G_m)$
under the null $F=G$ for independent samples with different sizes. Independently, looking for a feasible statistical 
way of relaxing the idea of ``treatment improvement" 
underlying stochastic dominance, \cite{Alv2017} introduced  (\ref{gamma}) as {a} na\"{\i}f index 
to measure deviation from stochastic dominance,  $F\leq_{st}G$ and provided some asymptotic 
theory for the empirical index, for the case of d.f.'s with a single crossing point (the typical case
in a location-scale family setting). In the same line, \cite{ZhuangEtAl2019} adapted the theory to cover even a finite number of crosses between the d.f.'s, under the additional assumption of  an exponential density ratio model and using semiparametric estimates of the quantile functions.

Here, in a wide setting, we provide a complete set of distributional limit results for Galton's rank order statistic, showing the  complex  panorama of the asymptotic behaviour of  $\gamma(F_n,G_m)$. In particular,  we  pursuit on the goal of analyzing the scarcely treated case of a finite number of contact points {between $F^{-1}$ and $G^{-1}$}, leading to a sound study of the local behaviour at every isolated contact point between quantile functions. This focuses on the consideration of the ``contact intensity" (to be properly defined), which exceeds the merely visual scope of crossing points of smooth enough curves and presents certain similarities with concepts lying in Stochastic Geometry. That contact intensity relies on the existence of a local {Lypschitzian} reparameterization of a curve in terms of the other.

 \subsection{Main results} 
Next we introduce the basic concepts we handle and explain the main results, whose proofs are deferred to Sections \ref{asymptoticresults} and \ref{puntos_separados}.

{Intuitively, the asymptotic behaviour of  $\gamma(F_n,G_m)$ depends on the size of the contact set, namely, the set
\begin{equation}\label{def_Gamma}
\Gamma:=\{t: F^{-1}(t)=G^{-1}(t)\}.
\end{equation}
For equal sized samples this was already observed in \cite{Gross}. We note that, since the index $\gamma$ is invariant with respect to strictly increasing transformations, the set $\Gamma$ could be equivalently expressed, in regular cases, as $\{t  : \varLambda(F^{-1}(t))=0\}$, where $\varLambda(x):=G^{-1}(F(x))-x$ is the shift function introduced in \cite{Doksum} as a richer alternative to the difference of means for comparing two continuous d.f.'s. The analysis of the  Q-Q process associated to $\varLambda$ was done in \cite{Aly}, under smoothness assumptions, through strong approximations. Yet, intuition may fail without some regularity conditions and, as we show in this work, $\Gamma$ is not really the right  set to look at. In fact, the asymptotic analysis of Galton's rank statistic is better handled in terms of the alternative shift function $h(t):=F_G(t)-t$, underlying the associated P-P process considered in \cite{Aly2}. Here, and throughout this work, we denote $F_G:=F\circ G^{-1}$ (similarly, $G_F=G\circ F^{-1}$) and
\begin{equation}\label{def_Gamma_tilde}
\tilde\Gamma:=\{t: F_G(t)=t\}.
\end{equation}
We observe that if $F$ and $G$ are continuous, then $\tilde\Gamma=\Gamma$. However, these sets can be quite different: for $F=G$,  a Bernoulli distribution  with mean $p$,  we have  $\Gamma=[0,1]$, while $\tilde\Gamma=\{1-p,1\}$.
By focusing on the `right' choice of contact set, our results go beyond the cases that could be treated  from the analyses in \cite{Aly} and \cite{Aly2}. In fact, we provide 
necessary and sufficient conditions for the a.s. consistency of $\gamma(F_n,G_m)$ \textit{without any smoothness assumption}:

\begin{Theo}\label{consistency} 
Let $F, G$ be arbitrary  d.f.'s.  Then $\gamma(F_n,G_m)\convs \gamma(F,G),  \mbox{ as } n,m \to \infty$ if and only if $\ell(\tilde\Gamma)=0$. 
\end{Theo}
A similar result holds for the one-sample statistic, $\gamma(F_n,G)$. 

As we see from Theorem \ref{consistency}, if $\ell(\tilde\Gamma)>0$, $\gamma(F_n,G_m)$ (or $\gamma(F_n,G)$) are not consistent estimators of $\gamma(F,G)$. In this case,  we provide a completely general result about the asymptotic behaviour of $\gamma(F_n,G)$: 
\begin{Theo}\label{one-side}
Let $F, G$ be arbitrary d.f.'s. Then $$\gamma(F_n,G)-\gamma(F,G) \convw \ell\left\{t\in \tilde\Gamma: B(t)>0\right\},$$ as $n\to \infty$, where $B$ is a standard Brownian bridge on $[0,1]$.
\end{Theo}

Still in the case $\ell(\tilde\Gamma)>0$, we prove weak convergence of the two-sample statistic, $\gamma(F_n,G_m)$, under mild assumptions.
This problem was also treated in \cite{Gross} for equal sample sizes ($n=m$), through combinatorial arguments and the method of moments.
That combinatorial approach seems to be inappropiate  to handle the case of unequal sample sizes. Additionally, our version yields a 
simple representation of the limit law.

\begin{Theo}\label{Theo_Noconsistency}
Let $F,G$ be  d.f.'s such that $F_G$ is Lipschitz. If $B$ is a standard Brownian bridge on $[0,1]$ and  $m,n\to \infty$ satisfy $0<\lim\inf\frac{n}{m+n}\leq\lim\sup\frac{n}{m+n}<1$, then
$$
\gamma(F_{n},G_{m})- \gamma(F,G)
\convw
\ell\{t \in \tilde\Gamma: B(t)>0\}.
$$
\end{Theo}

It should be noted that the limiting distribution in Theorems \ref{one-side} or \ref{Theo_Noconsistency} is non-degenerate if and only if $\ell(\tilde{\Gamma})>0$ (see Lemma \ref{nota_puente} in Section \ref{asymptoticresults}). If $\tilde{\Gamma}=(0,1)$, a celebrated result by Paul L\'evy (see Section 8, $2^o$ in \cite{Levy} or  p. 85-86 in \cite{Billingsley}) is that if $B$ is a standard Brownian bridge on $[0,1]$, then 
$$P\left(\ell\{t\in [0,1]: B(t)>0\}\leq x\right)=x, \mbox{ for every } x \in [0,1].$$
From this and Theorem \ref{Theo_Noconsistency} we recover, asymptotically, the classical result for the case $m=n$ and continuous $G=F$ (recall that in this case $\gamma(F_n,G_n)$ is uniformly distributed over $\{0,\frac 1{n},\ldots,\frac{n-1}n,1\}$; continuity of $F$ ensures that $F(F^{-1}(t))=t$ for every $t\in(0,1)$ and Theorem \ref{Theo_Noconsistency} applies with $\tilde{\Gamma}=(0,1)$).

When $\ell(\tilde{\Gamma})=0$ the limiting distribution in Theorems \ref{one-side} and \ref{Theo_Noconsistency} is degenerated at 0. In Section \ref{puntos_separados} we obtain non-degenerated limiting distributions, with different rates, when the contact set consists of a finite collection 
of  contact points. The key is the local asymptotic behaviour of $F_{n}^{-1}-G_{m}^{-1}$ around these influential points.
To avoid unnecessary smoothness assumptions here, we must consider contact points between nondecreasing functions in a generalized sense, including virtual contact points: those corresponding to contacts  between the vertical segments joining lateral limits at discontinuity points. Since quantile functions are left continuous, the following definition includes all these contact points.

\begin{Defi} \label{Defin.ContactPoint}
We say that $t\in (0,1)$ is a (generalized) contact point between $F^{-1}$ and $G^{-1}$ if either (i) $F^{-1}(t)=G^{-1}(t)$ or (ii) $F^{-1}(t)<G^{-1}(t)\leq F^{-1}(t+)$ or (iii) $G^{-1}(t)<F^{-1}(t)\leq G^{-1}(t+)$. The set of contact points will be denoted by $\Gamma^*.$
\end{Defi}

In Section \ref{Sec.Notation} we analyze the points in $\Gamma^*$ in detail. We show in particular (Proposition \ref{clave0}) that these generalized contact points are exactly the generalized contact points between the identity and the transforms $F\circ G^{-1}$ or $G\circ F^{-1}$. Notice that these non-decreasing functions can present also left-jump discontinuities which apparently would lead to additional virtual contact points between them and the identity. However, as we show in Proposition \ref{left_virtual}, for these functions ($F_G$ and $G_F$ vs. the identity), the same cases considered in Definition \ref{Defin.ContactPoint}, just including the contact points in the strict sense and those corresponding to right-jump discontinuities, suffice.

In some cases, and when it makes sense, quantile functions and the considered transforms can be extended by continuity to $0$ and $1$ and this, in turn, allows us to consider $0$ or $1$ as contact points (only in the strict sense). Hereafter we will distinguish \textit{extremal} contact points (0 or 1 when they are contact points) 
and \textit{inner} contact points (any other contact point).

If $\Gamma^*$ is finite, Corollary \ref{localizacion2} in Section \ref{puntos_separados} shows that the analysis of the asymptotic behaviour of $\gamma(F_n,G_m)$ boils down to the analysis of the {localized} measure of the set where the first sample quantile function exceeds  the other, namely, of
\begin{equation}
\label{Eq.Statistic.2}
\ell_{n,m}^{t_0} :=
\ell\left( \{ F_n^{-1}>G_m^{-1}\} \cap(t_0-\eta,t_0+\eta)\right)
-
\ell\left( \{F^{-1}>G^{-1}\}\cap(t_0-\eta,t_0+\eta)\right), 
\end{equation}
for $t_0\in \Gamma^*$, assuming that $\eta >0$ is small enough to ensure that $t_0$ is the only point in $\Gamma^*\cap(t_0-\eta,t_0+\eta)$ (as we will see, asymptotically, $\ell_{n,m}^{t_0}$ does not depend on $\eta$, hence, we do not include it in our  notation). We relate this behaviour to the character, position and intensity of the contact, in the following sense. For $t_0$ 
such that $F_G(t_0)=t_0$, set $\mathcal H:=\{h: t_0+h\in [0,1]\}$ and consider  the function $\Delta: \mathcal H \to \Rea$:
\begin{equation} \label{Eq.Delta}
\Delta (h) := F_G(t_0+h) - t_0 - h. 
\end{equation}
We will assume that {$\Delta$ is locally Lipschitz at $0$ (equivalently, that $F_G$ 
is locally Lipschitz at $t_0$) plus a higher order expansion, possibly on the positive and the negative  sides (for extremal contact points only one  expansion makes sense). More precisely, we will assume additionally {to the Lipschitz property} that} there exist $\eta>0$,  $r_L=r_L(t_0), r_R=r_R(t_0)\geq 1$  and $C_L=C_L(t_0)\ne 0, C_R=C_R(t_0)\ne 0$ such that
\begin{equation}
\label{Eq.DesarrolloDelta}
\Delta(h)=\left\{
\begin{matrix} 
C_L |h|^{r_L}  + o(|h|^{r_L}),& & \mbox{ if } h \in (-\eta,0),
\\[2mm]
C_R |h|^{r_R}  + o(|h|^{r_R}), & &\mbox{ if } h \in (0,\eta).
\end{matrix}
\right.
\end{equation}
In these cases, we will say that $r_L$ (resp. $r_R$) is the {\it  intensity or order of the contact on the left} (resp. {\it on the right}) between $F^{-1}$ and $G^{-1}$ at $t_0$. We observe that the assumptions imply that, for small enough $\eta$, $\sign(F_G(t)-t)=\sign(C_L)$ on $(t_0-\eta,t_0)$ and $\sign(F_G(t)-t)=\sign(C_R)$ on $(t_0,t_0+\eta)$. 
A point $t_0$ satisfying these conditions will be called a \textit{regular contact point}.

For integer $r_L$ and $r_R$, expression \eqref{Eq.DesarrolloDelta} is a kind of left- and right- Taylor expansion.
However, $r_L$ and $r_R$ are not necessarily integer in the definition above. 
We can classify regular contact points as {\it crossing points} (the case $\sign(C_L)\ne \sign(C_R)$) or {\it tangency points} 
(if $\sign(C_L)=\sign(C_R)$).
Notice also that under a proper Taylor expansion, $t_0$ is a crossing or tangency point depending only on whether $r_L=r_R$ is odd or even, while decomposition \eqref{Eq.DesarrolloDelta} allows to have a crossing point with odd $r_L=r_R$ or a tangency point with even  $r_L=r_R$.

We must stress that \eqref{Eq.DesarrolloDelta} does not necessarily imply smoothness conditions on $F^{-1}$ or $G^{-1}$. As an example, consider the case when $F^{-1}(t_0)<G^{-1}(t_0)\leq G^{-1}(t_0+)< F^{-1}(t_0+)$. Then $t_0$ is a discontinuity point of $F^{-1}$ (maybe also of $G^{-1}$), but $F_G$ is then locally constant, namely, $F_G(t)=t_0$ for $t$ close enough to $t_0$ and \eqref{Eq.DesarrolloDelta} holds with $r_L=r_R=1$, $C_R=-C_L=-1$. We should also note that, while \eqref{Eq.DesarrolloDelta} excludes discontinuity points for $F_G$, in particular, virtual contact points between $F_G$ and the identity, our approach allows to handle these points in a rather straigthforward way (see \eqref{virtual.crossing}, \eqref{tangencia.virtual} and Theorem \ref{finite.support}). Finally, we note that while \eqref{Eq.DesarrolloDelta} requires the contact orders to be at least 1, lower orders can also be considered. If, for instance, $\Delta(h)=\mbox{sgn}(h)|h|^r$, with $0<r<1$, then $\Delta$ is not Lipschitz around 0, but $G_F$ is and, under some additional assumptions, the local behaviour can be studied through $\tilde{\ell}_{m,n}^{t_0}$, the version of $\ell_{n,m}^{t_0}$ in which the roles of the $X$ and $Y$ samples are exchanged (see the comments before the proof of Theorem \ref{pieces}).

For 
a compact description of the limit distribution for the terms $\ell_{n,m}^{t_0}$ we consider independent random elements $B_1,B_2,W_0, W_1,\{\xi_{1,n}\}_{n\geq 1}, \{\xi_{2,n}\}_{n\geq 1}, \{\xi_{3,n}\}_{n\geq 1}$, $\{\xi_{4,n}\}_{n\geq 1}$,  where $B_i$ are Brownian bridges on $[0,1]$, $W_i$ are Brownian motions on $[0,\infty)$ and $\{\xi_{i,n}\}_{n\geq 1}$ sequences of i.i.d. exponential r.v.'s with unit mean. 
We set $S^i_{k}:=\xi_{i,1}+\cdots+\xi_{i,k}$, $k\geq 1, i=1,\ldots,4$. We fix $\lambda \in (0,1)$ and  set $B_\lambda:=\frac {1} {\sqrt{\lambda}}B_1  - \frac{1}{ \sqrt{1-\lambda}}B_2$. 

We consider $r_L,r_R\geq 1$ and denote $r_0:=\max(r_L,r_R)$. Also, for real numbers $a,b$,  we will use the notation $a^{{\mbox{\tiny sgn}(b)}}$ for  $a^+$ (the positive part of $a$)) either $a^-$ (the negative part) depending on whether $b>0$ or $b<0.$ For $t_0\in(0,1)$ we define 
\begin{eqnarray}\nonumber
T_{r_L,r_R}(t_0;C_L,C_R)&:=&\sign(C_L)\Big({\textstyle \frac{(B_\lambda(t_0))^{{\mbox{\tiny sgn}(C_L)}}}{|C_L|}}\Big)^{1/r_0}I(r_L=r_0)\\
\label{locallimit}
&&+\, 
\sign(C_R)\Big({\textstyle \frac{(B_\lambda(t_0))^{{\mbox{\tiny sgn}(C_R)}}}{|C_R|}}\Big)^{1/r_0}I(r_R=r_0),
\end{eqnarray}
when $r_0>1$ or $r_0=1$ and $C_RC_L>0$, while
\begin{eqnarray}\label{locallimit1}
T_{1,1}(t_0;C_L,C_R)&:=&{\textstyle \frac{(B_\lambda(t_0))^{{\mbox{\tiny sgn}(C_L)}}}{C_L}}+{\textstyle \frac{(B_\lambda(t_0))^{{\mbox{\tiny sgn}(C_R)}}}{C_R}}+\sign(C_L){\textstyle \frac{B_2(t_0)}{\sqrt{1-\lambda}}},
\end{eqnarray}
when $C_LC_R<0$. 
Additionally, for $r_0>1$ {and $t_0=0,1$} we define
\begin{eqnarray}
\label{locallimit0r}
T_{r_0,r_0}({t_0};C,C):=&\sign(C)\ell\big\{y\in (0,\infty):\,\sign(C) W_{t_0}(y)>(\lambda (1-\lambda))^{1/2}|C|y^{r_0}\big\},
\end{eqnarray}
while in the case $r_0=1$ we set
\begin{eqnarray}\nonumber
T_{1,1}(0;C,C):=&\sign(C)\lambda (1-\lambda) \int_0^\infty I_{\big\{\sign(C) {  \lambda } S^2_{\lceil (1-\lambda)  y \rceil} >\sign(C)(1-\lambda)(1+C) S^1_{\lceil \lambda   y \rceil}\big\} }dy
\\
\label{locallimit01}
T_{1,1}(1;C,C):=&\sign(C)\lambda (1-\lambda) \int_0^\infty I_{\big\{\sign(C) {  \lambda } S^4_{\lceil (1-\lambda)  y \rceil} >\sign(C)(1-\lambda)(1+C) S^3_{\lceil \lambda   y \rceil}\big\} }dy.
\end{eqnarray}
The double subindex and the double $C$ are redundant for these extremal contact points, but allows to keep a simple notation.

We are ready to present the results describing the asymptotic behaviour of $\ell_{n,m}^{t_0}$ for regular contact points. Theorem \ref{pieces} deals with innner contact points, while extremal contact points are considered in Theorem \ref{Theo_0_1.1} (in fact, Theorem \ref{pieces} remains valid for extremal contact points, but the limit distribution is Dirac's measure on 0 in that case).

\begin{Theo}\label{pieces}
Assume $t_0$ is a regular inner contact point with contact orders $r_L=r_L(t_0),r_R=r_L(t_0)\geq 1$ and 
constants $C_L=C_L(t_0),C_R=C_R(t_0)$. If $r=\max(r_L,r_R)$ and $n, m\to \infty$ with $\frac n {n+m} \to \lambda \in (0,1),$ then,  for every small enough $\eta>0$, 
$$(n+m)^{\frac 1 {2r}}\ell_{n,m}^{t_0}\convw T_{r_L,r_R}(t_0;C_L,C_R).$$
\end{Theo}

\begin{Theo}\label{Theo_0_1.1} 
Assume $t_0\in\{0,1\}$ is regular with contact order $r\geq 1$ and constant $C$. If $m, n \to \infty$ with $\frac{n}{n+m}\to\lambda \in (0,1)$, then,  for every small enough $\eta>0$
\begin{equation}\label{distr_limite}
{(n+m)^{\frac 1 {2r-1}}}\ell_{n,m}^{t_0}
\convw T_{r,r}(t_0;C,C).
\end{equation}
\end{Theo}

We see from Theorems \ref{pieces} and \ref{Theo_0_1.1} that, with the same contact intensities, $\ell_{n,m}^{t_0}$ vanishes faster for extremal contact points. In Subsection \ref{examplessec} we provide examples of extremal contact points for which $\ell_{n,m}^{t_0}$ converges at rate $(n+m)^{-c}$ for every $c\in(0,1]$ and of inner contact points for which the rate is $(n+m)^{-c}$, $c\in (0,\frac 1 2]$. Another distinctive feature of the limiting distributions for inner contact points is that for crossing points (those with $C_L(t_0)C_R(t_0)<0$) the limiting distribution takes positive and negative values with positive probabilities.  If $t_0$ is a tangency point ($C_L(t_0)C_R(t_0)>0$) then the limiting distribution is concentrated on $(0,\infty)$ or on $-(\infty,0)$.

The local asymptotic results in Theorems \ref{pieces} and \ref{Theo_0_1.1} can be strengthened to produce the following distributional limit
theorem for $\gamma(F_n,G_m)$.
\begin{Theo}\label{TCL.principal}
Assume $\Gamma^*=\{t_1,\ldots,t_k \}$ where $t_i$ is a regular contact point with intensities $r_L(t_i),r_R(t_i)$ and
constants $C_L{(t_i)}, C_R(t_i)$, $i=1,\ldots,k$. Set $r_{i}=\max(r_L(t_i),r_R(t_i))$ if $t_i\in (0,1)$, and $r_{i}=\max(r_L(t_i),r_R(t_i))-\frac 1 2$ if $t_i\in\{0,1\}$. Then, if $r_0=\max_{1\leq i\leq k} r_{i}$,
$$(n+m)^{\frac 1 {2r_0}}(\gamma(F_n,G_m)-\gamma(F,G)) \convw \sum_{i=1}^k I(r_{i}=r_0) T_{r_L(t_i),r_R(t_i)}(t_i;C_L(t_i),C_R(t_i)).$$
\end{Theo}

Theorem \ref{TCL.principal} shows that the rate of convergence of $\gamma(F_n,G_m)$ is determined by the maximal intensity of contact, and that only points with maximal intensity contribute to the limiting distribution, with adjustments to take into acount the different role of inner and extremal contact points. If there are extremal contact points then the rate of convergence can be $(n+m)^c$ for any $c\in(0,1]$. When there are only inner contact points the rate is $(n+m)^c$ with $c\in(0,\frac 1 2]$. The only case in which $\gamma(F_n,G_m)$ is asymptotically normal is when the inner contact points have intensity one, all of them with constants  $C_L=-C_R$ and there is no extremal contact point or its influence vanishes faster. 

\subsection{Organization of the paper}
The remaining sections of this work are organized as follows. Section \ref{Sec.Notation} includes  some  key results on quantile functions  and analyzes the structure of the contact sets. 
We will explicitly formulate several results on quantile functions. Some are classical, but, in fact it is not an easy task to find a comprehensive reference on quantile functions, with the notable exception of Appendix A in \cite{BobkovLedoux} on `Inverse Distribution Functions'. We observe that \cite{BobkovLedoux} is devoted to the analysis of convergence rates of Kantorovich transport distances between probability measures on the real line, which can be expressed in terms of quantile functions as $\int_0^1|F^{-1}(t)-G^{-1}(t)|^pdt,$ thus our problem corresponds to the limiting case $p=0$. Remarkably, this problem also encompasses a wide range of convergence rates.

In Section 3 we provide the proofs of  Theorems \ref{consistency}, \ref{one-side} and \ref{Theo_Noconsistency}.
Most of the limit theorems that we give for Galton's rank statistic are based on convenient representations of empirical quantile functions, 
combined with some type of strong approximation. 
Using representation \eqref{Eq.represent} below, we can derive limit theorems for Galton's rank statistic relying on strong approximations for  
uniform quantile processes, rather than using strong approximations for general quantile processes (as, for instance, in Chapter 6 in \cite{CsorgoHorvath}). This results in a significant gain in generality, since approximations for general quantile processes typically require strong smoothness assumptions (existence of densities plus additional conditions on them) that we can circumvent with this approach.

Section 4 gives the proofs of Theorems \ref{pieces}, \ref{Theo_0_1.1} and \ref{TCL.principal}. The key ingredients for this will be, as in Section 3, a convenient representation of the quantile processes and some application of strong approximations. With some simple localization results (Lemma \ref{localizacion} and Corollary \ref{localizacion2}) we see that the asymptotic behaviour of $\gamma(F_n,G_n)$ can be studied through that of the localized terms $\ell_{m,n}^{t_0}$ with $t_0$ in the contact set. Some results on the asymptotic independence between lower, central and upper order statistics allow then to complete the proof of Theorem \ref{TCL.principal}. Subsection \ref{examplessec} in that section provides  some examples of contact points with different positions and contact intensities. This subsection also includes a simplified version of Theorem \ref{TCL.principal}  under conditions that guarantee that $F_G$ is smooth (see Theorem \ref{dosmuestras_2}); and a further limit theorem (Theorem \ref{finite.support}) for the case  when $F$ and $G$ have  finite supports. This is an interesting example which can be handled with our approach even though the contact points here are not regular contact points. 

{
We include an Appendix with some additional material. The first part is devoted to some properties of the $F_G$ transform, including a technical discussion on conditions which guarantee that $F_G$ is Lipschitz or locally Lipschitz. Finally, we present a strong approximation result that we have used in several proofs.}

\subsection{Notation}

We end this Introduction with some words on notation. Through the paper $\mathcal L(X)$ will denote the law of the 
random vector or r.v. $X$.   We will consider a generic probability space $(\Omega,\sigma, P)$, where the involved random objects are defined. Given the (measurable) sets $A,B$, by $I_A$ we will denote the indicator function of  $A${ and $A\backslash B$ will denote the set $\{x\in A: \ x\notin B\}.$ 
} 
As before, $\ell$ will denote the Lebesgue measure on the  unit interval $(0,1)$. Convergences in the almost surely, in probability,  or in law (or weak) senses will respectively denoted by $\convs$, $\convp$, and $\convw$.  
Given a real value, $x$, we will  use $\lceil x\rceil$ to denote the smaller integer greater or equal than $x$, and $x^+:=\sup\{x,0\}$ and $x^-=-\inf\{x,0\}$. Also we use the notation  $f(x-):=\lim_{y\to x-}f(y)$ and $f(x+):=\lim_{y\to x+}f(y)$ for the lateral limits of a real function, $f$, whenever these limits exist, and $\sign(x)$  (defined, for real $x$, as $0$ if $x=0$ and $x/|x|$ otherwise). {Also recall that  for real numbers $a,b$, $a^{{\mbox{\tiny sgn}(b)}}$ will denote either $a^+$ or $a^-$ depending on whether $b>0$ or $b<0.$}

Throughout, $X_1,\dots,X_n$ and $Y_1,\dots,Y_m$ will be independent samples of i.i.d. r.v.'s such that $\mathcal L(X_i)$ and $\mathcal L(Y_i)$ have respective d.f.'s $F$ and $G$. As above,  $F_n$ and $G_m$  will denote the respective sample d.f.'s based on the $X's$ and $Y's$ samples.  Occasionally, we will use  the superscript $\omega$ in functions computed from the sample values $X_i(\omega),i=1,\dots,n$ or $Y_j(\omega),j=1,\dots,m$, (for instance, the empirical d.f. $F_n^\omega$ or the empirical quantile function $(F_n^\omega)^{-1}$). Without loss of generality we can (and often do) assume that the samples have been obtained from independent $U(0,1)$ samples $U_1,\dots,U_n$ and $V_1,\dots,V_m$ through the transformations $X_i=F^{-1}(U_i), Y_j=G^{-1}(V_j).$  From now on, we will denote the empirical quantile functions of these uniform samples by $\mathbb U_n$ and $\mathbb V_m$. We have the obvious relations $F_n^{-1}=F^{-1}(\mathbb U_n), G_m^{-1}=G^{-1}(\mathbb V_m).$ Writing $u_n$ and $v_m$ for the quantile processes based on the $U_i$'s and the $V_j$'s, respectively, ($u_n(t)=\sqrt n (\mathbb U_n(t)-t)$ and similarly for $v_m$) we note that

\begin{equation}\label{Eq.represent}
F_n^{-1}(t)=F^{-1}{\textstyle \Big(t+\frac{u_n(t)}{\sqrt{n}}\Big)} \ \mbox{ and } \ G_m^{-1}(t)=G^{-1}{\textstyle \Big(t+\frac{v_m(t)}{\sqrt{m}}\Big).}
\end{equation}

As already noted, the limiting behaviour of Galton's rank statistic is best described in terms of the contact set between the identity and the  function  
\begin{equation}\label{h-func}
{F_G(t):=F(G^{-1}(t)).}
\end{equation}

We note that, while the role of $F$ and $G$ is symmetric in the definitions of $\Gamma$ and $\Gamma^*$, this is not true in the case of $\tilde{\Gamma}$. For a more clear description of the relations among these sets we sometimes write $\tilde{\Gamma}_F=\tilde{\Gamma}$, $\tilde{\Gamma}_G=\{t\in(0,1): G_F(t)=t\}$ and 
$\tilde{\Gamma}^*_F$ (resp. $\tilde{\Gamma}^*_G$) for the set of {generalized} contact points  between $F_G$ (resp. $G_F$) and the identity (see (\ref{Defin.ContactPoint2})). Obviously,
$\tilde{\Gamma}_F\subset\tilde{\Gamma}^*_F$.

\section{Quantile transforms and contact sets} \label{Sec.Notation}

Quantile functions defined as in (\ref{defq}) provide a useful description of probabilities on the real line in terms of nondecreasing, left-continuous functions on $(0,1)$. In fact, every nondecreasing left-continuous real function, $H$, defined on $(0,1)$ is the unique  quantile function associated to just a unique d.f.: as a dual relation to  (\ref{defq}), such a function $H$ is the  quantile function associated to the d.f.
\begin{equation}\label{deff}
{F(x) =\sup\{t\in (0,1): H(t)\leq x\}.}
\end{equation}
As already noted, it will be convenient at some points to extend $F^{-1}$ to $0$ and $1$ in the obvious way (hence, $F^{-1}(0):=F^{-1}(0+)$ and  $F^{-1}(1):=F^{-1}(1-)$).

In this section we present some relevant facts on the relation between quantile functions and the composite functions {$F_G$ defined in \eqref{h-func}} without any smoothness assumption on the d.f.'s. We must begin by stressing the fact that, in general, we cannot guarantee even lateral continuity of $F_G$ (that would we only guaranteed for $F(G^{-1}(t)-)$ on the left and for $F(G^{-1}(t+))$ on the right).  On the other hand, from the well known relation for $t\in (0,1)$: $t\leq F(x) \iff F^{-1}(t) \leq x$,  it is easy to see the relations
\begin{equation}\label{h-func1}
F_G(t)=\max\{s\in [0,1]: F^{-1}(s)\leq G^{-1}(t)\}, \ \mbox{ for } t\in (0,1)
\end{equation}
{
\begin{equation}\label{compuesta}
F^{-1}(t)>G^{-1}(s) \iff  t>  F_G(s)  \ \mbox{ for } t,s\in (0,1).
\end{equation}

We note that $F_G(t)$ and $G_F(t)$ could be different even when $t \in \Gamma$. This possibility is naturally related to
the behaviour of the composition $F(F^{-1}(t))$. Clearly, $F(F^{-1}(t)-)\leq t \leq F_{F}(t)$,
thus $F_{F}(t)=t$ when $F$ is continuous  at $F^{-1}(t)$, but this could fail otherwise. More precisely, for $t\in (0,1)$, we have that   
$F_{F}(t)=t \ \mbox{ is equivalent to } \ t\in \mbox{Im}(F)$,  where $\mbox{Im}(F):=\{F(x) , x\in \Rea\}$ 
(see Lemma A.3 in \cite{BobkovLedoux}).  Now let $t\in (0,1)\cap \Gamma$: if  $t\in \mbox{Im}(F)$, then $t=F_F(t)=F_G(t)$, while if $t\notin \mbox{Im}(F)$, then $t<F_F(t)=F_G(t)$. We collect these facts  and some easy consequences {
for further reference} in the following  {
lemma}. 
 
\begin{Lemm}\label{Imagen}
Let $F,G$ be arbitrary d.f.'s. {
For $t\in (0,1)$, with the above notation, we have:
\begin{itemize}
\item[a)] If $t\in \Gamma \cap \mbox{Im}(F)$, then $F_G(t)=t$. 
\item[b)]  If  $t\in \Gamma\backslash\mbox{Im}(F),$ then $F_G(t)>t$.
\item[c)] If $F_G(t)=t$, then either $t\in \Gamma$ or $F^{-1}(t)<G^{-1}(t)$.
\item[d)] If $F_G(t)=t$ and $G_F(t)=t$, then $t\in \Gamma.$
\end{itemize}
}
\end{Lemm}
 
The conclusions in Lemma \ref{Imagen} can be rewritten with the notation of \eqref{def_Gamma} and \eqref{def_Gamma_tilde}. Item \textit{a)}, for instance, becomes $\Gamma\cap \mbox{Im}(F)\subset \tilde{\Gamma}$.
Generalized contact points in the sense of Definition \ref{Defin.ContactPoint} (that is, points in $\Gamma^*$) can also be characterized in terms of the composite functions $F_G$ and $G_F$.

As already noted in the Introduction, the consideration of virtual contact points associated to left-jump discontinuities of $F_G$ is not necessary because they are in fact contact points in the strict sense or they are associated to right-jump discontinuities of $G_F$ (see Proposition \ref{left_virtual}).  In consequence, we consider a point $t_0\in (0,1)$ as a contact point of $F_G$ and the identity whenever
\begin{equation}\label{Defin.ContactPoint2}
F_G(t_0)=t_0, \mbox{ or } F_G(t_0)<t_0 \leq F_G(t_0+)
\end{equation}

Note that the virtual contact condition $F_G(t_0)<t_0 \leq F_G(t_0+)$ is equivalent to $F_G(t_0)<t_0 \leq F_G(s)$ for all $s>t_0$, hence also to $G^{-1}(t_0)<F^{-1}(t_0)\leq G^{-1}(t_0+)$, which is condition \textit{(iii)} in Definition \ref{Defin.ContactPoint}. Therefore, we have shown that

\begin{Prop}\label{clave0}
The virtual contact points of $F^{-1}$ and $G^{-1}$ are exactly the virtual contact points of $F_G$ or $G_F$ with the identity.
\end{Prop}

 We note that $\tilde{\Gamma}^*_F\backslash\tilde{\Gamma}_F$ is contained in the set of discontinuity points of the nondecreasing function $F_G$, which must be at most countable. Hence,  $\ell(\tilde{\Gamma}^*_F\backslash\tilde{\Gamma}_F)=0$ and $\ell(\tilde{\Gamma}^*_F)=\ell(\tilde{\Gamma}_F)$. Proposition \ref{clave0} means that $\Gamma^*\backslash\Gamma=(\tilde{\Gamma}^*_F\backslash\tilde{\Gamma}_F)\cup(\tilde{\Gamma}^*_G\backslash\tilde{\Gamma}_G)$.
We explore next the situation for contact points in the strict sense.
\begin{Prop}\label{clave1}
If $t_0\in \Gamma\cap (0,1)$ then $F_G(t_0)=t_0,$ or $G_F(t_0)=t_0$ (that is $t_0 \in \tilde\Gamma_F\cup \tilde\Gamma_G$), or the set $\{t\in (0,1): F^{-1}(t)=G^{-1}(t)=F^{-1}(t_0)\}$ is a non-degenerate interval (hence, in the latter case,  the point $x_0=F^{-1}(t_0)$ is a common discontinuity point of $F$ and $G$ and $t_0$ cannot be an isolated element of $\Gamma^*$).
\end{Prop} 
\textbf{Proof.}  
It is easy to see that if $t_0\in (0,1)\backslash(\mbox{Im}(F)\cup \mbox{Im}(G))$ satisfies $F^{-1}(t_0)=G^{-1}(t_0)$, then the point $x_0=F^{-1}(t_0)$ would have positive mass under both distributions, hence the set $\{t\in (0,1): F^{-1}(t)=G^{-1}(t)=F^{-1}(t_0)\}$ is a non-degenerate interval.  Any other point in $\Gamma\cap (0,1)$ must belong to   $\mbox{Im}(F)\cup \mbox{Im}(G)$ and, by Lemma \ref{Imagen},  must satisfy either $F_G(t_0)=t_0$ or $G_F(t_0)=t_0.$ 
\FIN 
   
\begin{Prop}\label{clave2}
Let $t_0\in (0,1)$ be such that $t_0\in \tilde\Gamma_F\cup \tilde\Gamma_G$, that is $F_G(t_0)=t_0$ or  $G_F(t_0)=t_0$. Then $t_0 \in \Gamma^*$ ($t_0$ is a contact point between $F^{-1}$ and $G^{-1}$).\end{Prop} \textbf{Proof.}  
For any $t_0\in(0,1)$ such that  $F_G(t_0)=t_0$ (the case $G_F(t_0)=t_0$ is identical),  we must have one of the following exclusive  possibilities:
\begin{itemize}
\item[i)]  $G_F(t_0)<t_0$, and then $F^{-1}(t_0)<G^{-1}(t_0)$, and
\item[ii)] $F_G(t_0)= t_0= G_F(t_0)$, or $F_G(t_0)=t_0<G_F(t_0)$,  which lead to $F^{-1}(t_0)=G^{-1}(t_0).$
\end{itemize}
If i) holds, then we would have $G^{-1}(t_0)\leq F^{-1}(t_0+)$ (this follows easily from the fact that the strict inequality $G^{-1}(t_0)> F^{-1}(t_0+)$ would imply $F_G(t_0)=F(G^{-1}(t_0))>t_0$). Hence,  i) implies $F^{-1}(t_0)<G^{-1}(t_0)\leq F^{-1}(t_0+)$. 
  \FIN
  
 The next proposition  shows that it is not necessary to consider contact points associated to left-discontinuities.
   
  \begin{Prop}\label{left_virtual}
Let $t_0 \in (0,1).$ If $F_G(t_0-)\leq t_0<F_G(t_0),$ then $G_F(t_0)\leq t_0 \leq G_F(t_0+),$ or  the point $x_0=G^{-1}(t_0)$ is a common discontinuity point of $F$ and $G$ and $t_0$ cannot be an isolated element of $\Gamma^*$.
\end{Prop}
\textbf{Proof.}  
If we suppose that $G^{-1}(t)=G^{-1}(t_0)$ for some $t<t_0,$ then for every sequence $\{t_n\}$ such that $t_n\to t_0-,$ $F(G^{-1}(t_n))=F(G^{-1}(t_0))$ will hold eventually, thus leading to the absurd $F_G(t_0-)=F_G(t_0).$ Therefore it must be $x_0:=G^{-1}(t_0)>G^{-1}(t)$ for every $t<t_0,$ and $F_G(t_0-)=F(G^{-1}(t_0)-).$ Moreover, the discontinuity of $F$ and its link with $F^{-1}$ easily show that  $F^{-1}(t)=x_0 \mbox{ if } t\in (F(x_0-),F(x_0)]$. Now, on the first hand, from the hypothesis we obtain $F^{-1}(t_0)\leq x_0$ and $F^{-1}(s)=x_0$ for every $s\in (t_0,F(x_0))$, hence 
\begin{equation}\label{otramas}
G_F(t_0+)=G(x_0)=G(G^{-1}(t_0))\geq t_0.
\end{equation}
On the other hand, from the relation $t_0<F_G(t_0)$ we obtain $F^{-1}(t_0)\leq G^{-1}(t_0)$, hence   $F^{-1}(t_0)=G^{-1}(t_0)$ or, alternatively, $F^{-1}(t_0)<G^{-1}(t_0)$  what gives $G_F(t_0)< t_0$. This relation and (\ref{otramas}) imply that  $G_F(t_0)< t_0 \leq G_F(t_0+)$. 
Finally, if $F^{-1}(t_0)=G^{-1}(t_0)$ and $x_0=F^{-1}(t_0)$ is a continuity point of $G$, from (\ref{otramas}) we obtain that $G_F(t_0)= t_0$ what proves the result.
\FIN

{ We conclude this section with some easy consequences of the last results. }
\begin{Coro} \label{Coro.EquivContactos}
Let $t_0\in (0,1)$ such that $F_G(t_0)=t_0$ (resp. $G_F(t_0)=t_0$) then $t_0$ is a contact point (possibly virtual) between $G_F$ (resp. $F_G$) and the identity.
\end{Coro}
Corollary \ref{Coro.EquivContactos} states  that $\tilde{\Gamma}_F\subset \tilde{\Gamma}^*_G$. From the comments after Proposition \ref{clave0} we see that
$\ell(\tilde{\Gamma}_F)\leq \ell(\tilde{\Gamma}^*_G)=\ell(\tilde{\Gamma}_G)$. The same argument shows that $\ell(\tilde{\Gamma}_G)\leq\ell(\tilde{\Gamma}_F)$, hence $\ell(\tilde{\Gamma}_F)= \ell(\tilde{\Gamma}_G)$. This means, in particular,  that the roles of $F$ and $G$ in the condition $\ell(\tilde{\Gamma})=0$ in Theorems \ref{consistency}, \ref{one-side} and \ref{Theo_Noconsistency} are completely symmetric.

{ \begin{Prop}\label{Prop.clave012}
If $\Gamma^*$ is finite then $\Gamma^*=\tilde{\Gamma}^*_F\cup \tilde{\Gamma}^*_G$. In particular, $\tilde{\Gamma}^*_F$ and $\tilde{\Gamma}^*_G$ are finite.
\end{Prop}

We remark that, while $\tilde{\Gamma}^*_F\cup \tilde{\Gamma}^*_G\subset \Gamma^*$ always holds (this follows from Proposition \ref{clave2}), the set $\Gamma^*$ can be much bigger that $\tilde{\Gamma}^*_F\cup \tilde{\Gamma}^*_G$ (recall the comments in the Introduction; the case $G=F$, with $F$ the d.f. of the Bernoulli law with mean $p$ gives a simple example of this).

In Section \ref{puntos_separados} we prove distributional limit theorems for $\gamma(F_n,G_m)$ under the assumption that $\Gamma^*$ is finite, say,
$\Gamma^*=\{t_1<\cdots<t_r \}$. The differences $F^{-1}(t)-G^{-1}(t)$ must have constant sign in the open intervals $(t_i,t_{i+1})$   (the same happens in $(0,t_1)$ or
$(t_r,1)$ if $0$ or $1$ are not contact points). The next result will enable us to focus on neighbourhoods of isolated contact points to study $\gamma(F_n,G_m)$.

\begin{Lemm} \label{Lemm.SignoConstante}
Assume $0<a\leq b <1$ are such that $[a,b]\cap \Gamma{
^*} = \emptyset$, and also that $\mbox{sgn}(F^{-1}(t)-G^{-1}(t))>0$ (resp. $\mbox{sgn}(F^{-1}(t)-G^{-1})(t)<0$) for every $t\in[a,b]$. Then there exists $\delta >0$ such that $F^{-1}(t)-G^{-1}(t+)>\delta$ (resp. $F^{-1}(t)-G^{-1}(t+) < - \delta$)  for every $t \in [a,b]$.
\end{Lemm}
\noindent \textbf{Proof}: Let us consider the case $\mbox{sgn}(F^{-1}(t)-G^{-1}(t))>0$.
Assume, on the contrary, that there exist a sequence $\{t_k\}\subset [a,b]$ such that $F^{-1}(t_k)-G^{-1}(t_k+) \to 0$. In this case it is possible to choose $\{t^*_k\}\subset (0,1)$ such that, $t_k<t^*_k$, $t_k-t^*_k\to 0$,   $F^{-1}(t_k+)-F^{-1}(t_k^*) \to 0$ and  $G^{-1}(t_k+)-G^{-1}(t_k^*) \to 0$. Since $[a,b]$ is compact, we can assume that $\{t_k\}$ converges. Then also $\{t^*_k\}$ converges. We write $t_0 \in [a,b]$ for the common limit. 
By taking subsequences, if necessary, we can also assume that both sequences are monotone.
 
Now, we only need to consider   four possible cases. If, for instance,  $\{t_k\}$ and $\{t_k^*\}$ are increasing, then we would obtain that  $F^{-1}(t_0) = G^{-1}(t_0)$ which is impossible by assumption. If the sequence $\{t_k\}$  is increasing and $\{t_k^*\}$ is decreasing,  then $F^{-1}(t_k)\to F^{-1}(t_0)$ and $G^{-1}(t_k^*)\to G^{-1}(t_0+)$, and we would have that $F^{-1}(t_0) = G^{-1}(t_0+)$ and, consequently, $t_0$ would be a contact point what is not possible either, because $[a,b]\cap \Gamma^* = \emptyset$. The two remaining cases lead to similar contradictions.
  \FIN

We conclude this section with two observations. First, we note that $\sign(t-F_G(t))=\sign(F^{-1}(t)-G^{-1}(t))$ for every $t\notin \Gamma^*$. To check this  recall relation (\ref{compuesta}),  giving that {$F^{-1}(t)>G^{-1}(t)$ if and only if $t>F_G(t)$. This also implies that} $F^{-1}(t)\leq G^{-1}(t)$ if and only if $t\leq F_G(t)$, but $t=F_G(t)$ cannot happen if $t\notin\Gamma^*$ (Proposition \ref{clave2}). On the other hand, if $t< F_G(t)$ then
$F^{-1}(t)\leq G^{-1}(t)$ but, again, $F^{-1}(t)=G^{-1}(t)$ is not possible if $x\notin \Gamma^*$. This means that $\sign(t-F_G(t))$ is constant in the intervals $(t_i,t_{i+1})$ as above. 
  
Our second observation arises from  the fact that every nondecreasing left-continuous real function, $H$, defined on $(0,1)$ is the quantile function associated to the d.f. given by (\ref{deff}). We can apply Lemma \ref{Lemm.SignoConstante} to the quantile function $H(t)=F_G(t-)$ and the identity and conclude, for instance, that in a compact interval where $t-F_G(t)>0$ there exists some $\delta >0$ such that $t-F_G(t)\geq \delta$. We will exploit these facts in later sections.

}
}
 



\section{Consistency of Galton's rank order statistic}\label{asymptoticresults}

In this section we provide proofs of Theorems \ref{consistency}, \ref{one-side} and \ref{Theo_Noconsistency}. These results show that Galton's rank order statistic is a consistent estimator of the index $\gamma(F,G)=\ell\{t: F^{-1}(t)>G^{-1}(t)\}=
\ell\{t: t>F_G(t)\}$ if and only if the contact set $\tilde{\Gamma}=\{t: t=F_G(t)\}$ has zero Lebesgue measure. The key to the proof of Theorem \ref{consistency} is the following lemma. Here we denote $\Gamma_0:=\mbox{Im}(F)\cap \Gamma\cap (0,1)$.

\begin{Lemm}\label{LemmPrincipal}
Let $F, G$ be arbitrary  d.f.'s. With the notation  above, we have:
\begin{equation}\label{congamma}
\gamma(F_n,G_m)-\gamma(F,G)- \ell(\{F_n^{-1}>G_m^{-1}\}\cap \Gamma) \convs 0 \ \mbox{ as } n,m \to \infty,
\end{equation}
\begin{equation}\label{gamma0}
\gamma(F_n,G_m)-\gamma(F,G)- \ell(\{F_n^{-1}>G_m^{-1}\}\cap \Gamma_0) \convs 0 \ \mbox{ as } n,m \to \infty,
\end{equation}
and
\begin{equation}\label{fsubg}
\gamma(F_n,G_m)-\gamma(F,G)- \ell(\{F_n^{-1}>G_m^{-1}\}\cap \tilde\Gamma ) \convs 0 \ \mbox{ as } n,m \to \infty.
\end{equation}
\end{Lemm} 
 {\bf Proof.}  By right continuity, if $t \notin \mbox{Im}(F)$, then there exists $\delta_t>0$ such that $[t,t+\delta_t) \cap \mbox{Im}(F)=\emptyset$ and $F_F(s)=t+\delta_t$, for every $s \in [t,t+\delta_t]$. From this, it easy to see that there exists an at most countable family of disjoint intervals $I_k=[a_k,b_k)$, with $a_k <b_k$ which is a partition of the complement of $ \mbox{Im}(F)$ and $F_F(s)=b_k$,  for every $s \in [a_k,b_k)$.
 
  The Glivenko-Cantelli Theorem gives that for some $ \Omega_0 \in \sigma$, with $P( \Omega_0) =1$, if $\omega \in \Omega_0$, then
\[
\sup_t |F_n^\omega (t) - F(t)| \conv 0 \mbox{ and } \sup_t |G_m^\omega (t) - G(t)| \conv 0.
\]

Now, recalling the elementary Skorohod theorem  (see e.g. Lemma A.5 in \cite{BobkovLedoux}), for every  $\omega \in \Omega_0$, the set
\[
{\cal T}^\omega:=\left\{t\in (0,1) : \left(F_n^{\omega}\right)^{-1}(t) \to F^{-1}(t) \mbox{ and } \left(G_m^{\omega}\right)^{-1}(t) \to G^{-1}(t)\right\}\backslash \{a_1, a_2,\ldots\}
\]
 has Lebesgue measure one. Therefore,  if $\omega \in \Omega_0$,
\begin{eqnarray*}
\lefteqn{\gamma(F_n^\omega,G_m^\omega)-\gamma(F,G)- \ell(\{(F_n^\omega)^{-1}>(G_m^\omega)^{-1}\}\cap \Gamma )
}\\
&=&\ell \left[\{(F_n^\omega)^{-1}>(G_m^\omega)^{-1},F^{-1}< G^{-1}\}\cap {\cal T}^\omega \right]
-\ell \left[\{(F_n^\omega)^{-1}\leq(G_m^\omega)^{-1},F^{-1}> G^{-1}\}\cap {\cal T}^\omega)\right],
\end{eqnarray*}
which converges to 0 because both sets within brackets converge to the empty set. This proves (\ref{congamma}). 
To prove (\ref{gamma0}) we show that if $\omega \in \Omega_0$, then
\begin{equation} \label{Eq.Limitecs_1}
d_n:=\ell\left(\left\{(F_n^\omega)^{-1}>(G_m^\omega)^{-1}\right\}\cap \Gamma \right) - \ell\left(\left\{(F_n^\omega)^{-1}>(G_m^\omega)^{-1}\right\}\cap \Gamma _0\right) \to 0.
\end{equation}
To check this, notice that 
\begin{eqnarray*}
d_n &= &
\ell\left(\big\{(F_n^\omega)^{-1}>(G_m^\omega)^{-1}\big\}\cap {\cal T}^\omega \cap \Gamma \cap \big(\cup_k(a_k,b_k)\big) \right)
\\
&= &
\ell\left(\big\{t>(F_n^\omega)_{G_m^\omega}(t)\big\}\cap {\cal T}^\omega \cap \Gamma \cap \big(\cup_k(a_k,b_k)\big) \right). 
\end{eqnarray*}
Now, Glivenko-Cantelli again, and the construction of ${\mathcal T}^\omega$ yield that if $\omega \in \Omega_0$ and $t \in {\cal T}^\omega \cap \Gamma$, then
\begin{equation} \label{Eq.Limitecs_2}
0 =\lim_n \left| F_n^\omega\left[\left(G_m^\omega\right) ^{-1}(t)\right] -  F\left[\left(G_m^\omega\right) ^{-1}(t)\right]\right| \mbox{ and } \lim_n \left(G_m^\omega\right)^{-1}(t) = G^{-1}(t)= F^{-1}(t).
\end{equation}
From here, if $t \in {\cal T}^\omega \cap \Gamma\cap (a_k,b_k)$ for some $k$, then, eventually,
$
F[(G_m^\omega)^{-1}(t)] = b_k >t
$
which, combined with the first statement in \eqref{Eq.Limitecs_2}, makes eventually impossible that $t>(F_n^\omega)_{G_m^\omega}(t)$ and shows \eqref{Eq.Limitecs_1}.

The proof of (\ref{fsubg}) is now obvious taking into account that, from Lemma \ref{Imagen}
\begin{equation}\label{relations}
\Gamma_0\subset \tilde \Gamma \subset \Gamma \cup \{F^{-1}<G^{-1}\}.
\end{equation}
\FIN

\medskip
\noindent
{\bf Proof of Theorem \ref{consistency}.}  Sufficiency is a trivial consequence of Lemma \ref{LemmPrincipal}. To prove necessity, if $\gamma(F_n,G_m)\convs \gamma(F,G),  \mbox{ as } n,m \to \infty$, according to Lemma  \ref{LemmPrincipal}, we have that 
  \begin{equation}\label{Consistency0}
D_n:= \ell\left(\left\{F_n^{-1}>G_m^{-1}\right\}\cap \tilde\Gamma\right) \convs 0,
 \end{equation}
From (\ref{importante}), we have
\begin{equation*}
D_n=\ell\left(\left\{\mathbb U_n>F_G(\mathbb V_m)\right\}\cap \tilde\Gamma\right)                                               \geq  \ell\left(\left\{t: \mathbb U_n(t)> t\right\}\cap \left\{t: t \geq F_G(\mathbb V_m(t))\right\}\cap \tilde\Gamma \right) .
\end{equation*}
Now Fubini's theorem and independence between samples yield
\begin{equation}\label{Consistency2}
E [D_n] \geq \int_{\tilde\Gamma} P\left[ \mathbb U_n(t)> t\right] P\left[  t\geq F_G(\mathbb V_m(t))\right] dt \geq \int_{\tilde\Gamma} P\left[ \mathbb U_n(t)> t\right] P\left[  t> \mathbb V_m(t)\right] dt,
\end{equation}
where the last inequality follows from the fact that, since $F_G$ is nondecreasing, $F_G(t)=t$ and $\mathbb V_m(t)<t$ imply that $F_G(\mathbb V_m(t))\leq t.$ On the other hand, for every $t\in (0,1)$, the factors inside the integral converge to $1/2$. But, since $|D_n| \leq 1$,   \eqref{Consistency0} implies $E[D_n]\to 0$. This and convergence to 1/4 of the last integrand in \eqref{Consistency2} imply that $\ell (\tilde\Gamma)=0$.   
\FIN
}

Next, we give a proof of Theorem \ref{one-side}. We remark that our approach allows to handle this one-sample statistic without any smoothness assumption on $F$ or $G$. 

\medskip
\noindent
{\bf Proof of Theorem \ref{one-side}.}
Assuming, w.l.o.g., the construction in Theorem \ref{aproxBridge} we have
\begin{eqnarray*}
\hat \gamma_n&: =&\gamma(F_n,G)=\ell\left\{ t:F^{-1}\left(t+{\textstyle \frac{u_n(t)}{\sqrt n}}\right)>G^{-1}(t)\right\}  \\
&=&\ell\left\{t:t+{\textstyle \frac{u_n(t)}{\sqrt n}}>F_G(t)\right\}=\ell\left\{t:u_n(t)>{\sqrt n}(F_G(t)-t)\right\},
\end{eqnarray*}
and similarly 
$\gamma:= \gamma(F,G)=\ell\{t: F_G(t)-t<0\}$.
Therefore, we see that
\begin{equation*}
\hat\gamma_n-\gamma=\ell\left\{t: u_n(t)>{\sqrt n}(F_G(t)-t)\geq 0\right\}-\ell\left\{t: 0> \sqrt n(F_G(t)-t)\geq u_n(t)\right\}.
\end{equation*}

Obviously, for the Brownian bridges $B_n^F(t)$, 
\begin{eqnarray*}
\ell\left\{t: u_n(t)>{\sqrt n}(F_G(t)-t){>} 0\right\}&\leq& \ell\left\{t: B_n^F(t)+K{\textstyle \frac{\log n}{\sqrt n}}\geq \sqrt n(F_G(t)-t)>0\right\}
\\
&&+ \ell \left\{t: |B_n^F(t)-u_n(t)|>K{\textstyle \frac{\log n}{\sqrt n}}\right\}.
\end{eqnarray*}

By Theorem \ref{aproxBridge}, the last summand  eventually vanishes. For a fixed Brownian bridge $B(t)$ and $t\in (0,1)$ such that $F_G(t)-t>0$, we have $B(t)+K\frac{\log n}{\sqrt n}<\sqrt n(F_G(t)-t)$ eventually. This and the bounded convergence theorem imply that 
$$
\ell\left\{t: B(t)+K{\textstyle \frac{\log n}{\sqrt n}}\geq \sqrt n(F_G(t)-t)>0\right\}\convs 0.$$ 

As a result we obtain that $$\ell\left\{t: u_n(t)>{\sqrt n}(F_G(t)-t)>0\right\}\convp 0.$$

Similarly we see that $$\ell\left\{t: 0> \sqrt n(F_G(t)-t)\geq u_n(t)\right\} \convp 0$$
and conclude that 
\begin{equation}\label{one-side4}
\hat\gamma_n-\gamma=\ell\{t: u_n(t)\geq 0, F_G(t)=t\}+o_P(1).
\end{equation}

Next, we observe that, eventually,
\begin{eqnarray*}
\ell\left\{t \in \tilde\Gamma:B_n^F(t)-K{\textstyle \frac{\log n}{\sqrt n}}\geq 0\right\}&\leq& \ell\{t \in \tilde\Gamma: u_n(t)\geq 0\}
\\ 
&\leq& \ell \left\{t \in \tilde\Gamma: B_n^F(t)+K{\textstyle \frac{\log n}{\sqrt n}}\geq 0\right\}.
\end{eqnarray*}

Now, $$\ell \left\{t \in \tilde\Gamma: B(t)+K{\textstyle \frac{\log n}{\sqrt n}}\geq 0\right\}\conv \ell \left(t \in \tilde\Gamma: B(t)\geq 0\right),$$
and 
$$\ell \left\{t \in \tilde\Gamma: B(t)-K{\textstyle \frac{\log n}{\sqrt n}}\geq 0\right\}\conv \ell \left\{t \in \tilde\Gamma: B(t)\geq 0\right\}.$$
This and (\ref{one-side4}) show the announced result.
\FIN

We recall that the set involved in the limit law in the last result is $\tilde\Gamma$, which generally does not coincide with $\Gamma$ (see Lemma \ref{Imagen} and (\ref{relations}) for more details). For a better understanding of the links between Theorem \ref{consistency} and \ref{one-side}, we note that degeneracy in the limit law is equivalent to $\ell (\tilde\Gamma)=0$. This is an obvious consequence of the next, simple result.

\begin{Lemm}\label{nota_puente}
If $B(t)$ is a standard Brownian bridge on $[0,1]$, for any Borel set $A$ in $[0,1]$, the r.v. $\ell(\{B>0\}\cap A)$ is a.s. constant if and only if $\ell(A)=0$.
\end{Lemm}
{\bf Proof.}
{\rm
If $\ell(A)=0$ then, obviously, $\ell(\{B>0\}\cap A)=0$. Assume now that $\ell(A)>0$. It is well known that $\ell\{t\in [0,1]: B(t)=0\}=0$ (this follows easily from  Fubini's Theorem). Moreover, if $B$ is a Brownian bridge then $B=_d-B$.  
Hence, $\ell( \{B<0\}\cap A)=_d\ell( \{B>0\}\cap A)$, while   $\ell( \{B<0\}\cap A) +\ell( \{B>0\}\cap A)=\ell(A)$. This implies that 
$E(\ell( \{B>0\}\cap A))=\ell(A)/2$. Thus, if $\ell( \{B>0\}\cap A)$ were a.s. constant, that constant should equal $\ell(A)/2$. However, 
$\ell\{B>0\}$ stochastically dominates $\ell\{\{B>0\}\cap A\}$, and degeneracy on the value $\ell(A)/2$ would lead to the conclusion that the $U(0,1)$ law stochastically dominates Dirac's measure on $\ell(A)/2$, which cannot hold if $\ell(A)>0$.
\FIN

To deal with Galton's rank statistic in the two-sample case we must adapt the argument in the proof of Theorem \ref{one-side}.
This is done with Lemma \ref{Lem.BaseDeLoDeTasio}, which will play an important role in our development.  It relies on the strong approximation given in Theorem \ref{aproxBridge} in the Appendix. Given two real functions $f$ and $g$ and versions of independent sequences of Brownian bridges $\{B_{n}^F\}$,  $\{B_{m}^G\}$  and of uniform quantile processes, $u_n$ and $v_m$, as in  Theorem \ref{aproxBridge},  we set 
\begin{equation}\label{Eq.Las Tildes}
\begin{array}{rcl}
{f}_n(t):=f(t+\frac{u_n(t)}{\sqrt{n}}) &\mbox{ and  }& {g}_m := g(t+\frac{v_m(t)}{\sqrt{m}}),
\\
[2mm]
\tilde{f}_n(t):=f({\textstyle t +\frac{B_{n}^F(t)}{\sqrt{n}} }) &\mbox{ and  }& \tilde{g}_m := g({\textstyle t +\frac{B_{m}^G(t)}{\sqrt{m}} }).
\end{array}
\end{equation}

\begin{Lemm} \label{Lem.BaseDeLoDeTasio}
Consider $A\subset (0,1)$ such that $\ell(A)>0$. With the notation and construction of Theorem \ref{aproxBridge}, if we assume that $f,g$ are two real Lipschitz functions, then there exists $L>0$ such that, if $C_{n,m}:= L ({\textstyle \frac {\log n} n + \frac {\log m} m})$,  then whenever $n,m \to \infty$, eventually,
\begin{eqnarray}
\nonumber
\ell\big\{t\in A:\, \tilde{f}_n(t)>\tilde{g}_m(t)+C_{n,m}\big\}
&\leq &\ell\big\{t\in A:\, {f}_n(t)>{g}_m(t)\big\}
\\
\label{Eq.TgEnCero_2}
&\leq&  \ell\big\{t\in A:\, \tilde{f}_n(t) > \tilde{g}_m(t)-C_{n,m}\big\}.
\end{eqnarray}
\end{Lemm}

\noindent \textbf{Proof:} 
Since $f$ is Lipschitz, for $t \in A$ we have that
\begin{eqnarray*}
\big|f_n{\textstyle( t )}- \tilde{f}_n(t) \big|&=&
\big|f({\textstyle t +\frac{u_n(t)}{\sqrt{n}} }) -f({\textstyle t +\frac{B_{n}^F(t)}{\sqrt{n}} }) \big|
\leq  
\| f\|_{\mbox{\scriptsize Lip}} {\textstyle \frac{\|u_n-B_{n}^F\|_\infty}{\sqrt{n}}},
\end{eqnarray*}
with a similar bound for $\big|g_m{\textstyle( t )}- \tilde{g}_m(t) \big|$. These  bounds and (\ref{eventualbound2}) imply that on a probability one set, eventually,
\[
\sup_{t \in A} \left|(f_n-g_m)-(\tilde{f}_n-\tilde{g}_m)\right| \leq L ({\textstyle \frac {\log n} n + \frac {\log m} m}) = C_{n,m}
\]
for some positive constant $L$ (depending only on $f$ and $g$).
Observe that
{
\begin{eqnarray*}
\lefteqn{\ell\big\{t\in  A:\, f_n(t)>g_m(t)\big\} 
\ \leq  \ \ell\big\{t\in  A:\, \tilde{f}_n(t)>\tilde{g}_m(t)-C_{n,m}\big\}
}\hspace*{3cm}
\\
& & +\ell\big\{t\in  A:\, |(f_n(t)-\tilde{f}_n(t))-(g_{m}^{-1}(t)-\tilde{g}_m(t))|>C_{n,m}\big\},
\\
[2mm]
\lefteqn{
\ell\big\{t\in  A:\, \tilde{f}_n(t)>\tilde{g}_m(t)+C_{n,m}\big\} \ \leq \  \ell\big\{t\in A:\, {f}_n(t)>{g}_m(t)\big\}}\hspace*{3cm}\\
&&+\ell\big\{t\in  A:\, |(f_n(t)-\tilde{f}_n(t))-(g_{m}^{-1}(t)-\tilde{g}_m(t))|> C_{n,m}\big\}.
\end{eqnarray*}
On a probability one set the second summands on the last two upper bounds eventually vanish. Hence, on that probability one set,  (\ref{Eq.TgEnCero_2}) eventually holds.
\FIN

We will apply Lemma \ref{Lem.BaseDeLoDeTasio} to the cases in which $f=F^{-1}$ and $g=G^{-1}$ and when $f$ is the  identity and $g=F_G$ (see Section \ref{lips} in the Appendix for the analysis of the Lipschitz condition on $F_G$).

{
We end the section with the proof of  the two-sample analogue of Theorem \ref{one-side}. 
\

\medskip
\noindent
{\bf Proof of Theorem \ref{Theo_Noconsistency}.} By taking subsequences we can assume $\frac n {n+m}\to \lambda\in (0,1)$. Also, after Lemma \ref{LemmPrincipal}, it suffices to prove that 
$$
\ell\{t \in \tilde\Gamma: F_n^{-1}(t)>G_m^{-1}(t) \}\convw \ell \{t \in \tilde\Gamma: B(t)>0\},
$$
which, using the approximation in Theorem \ref{aproxBridge}  and Lemma \ref{Lem.BaseDeLoDeTasio}, will hold if 
\begin{equation}\label{primera}
\ell\big\{t\in \tilde\Gamma,{\textstyle t +\frac{B_{n}^F(t)}{\sqrt{n}} }>F_G({\textstyle t +\frac{B_{m}^G(t)}{\sqrt{m}} })+C_{n,m}\big\} \convw \ell \{t \in \tilde\Gamma: B(t)>0\}
\end{equation}
and
\begin{equation}\label{segunda}
\ell\big\{t \in \tilde\Gamma:\,  {\textstyle t +\frac{B_{n}^F(t)}{\sqrt{n}} }>F_G({\textstyle t +\frac{B_{m}^G(t)}{\sqrt{m}} })-C_{n,m}\big\} \convw \ell \{t \in \tilde\Gamma: B(t)>0\}.
\end{equation}
Both terms can be handled similarly, hence we will address here only (\ref{primera}). First, we note that $\ell\{t \in \tilde\Gamma: F_G(t)\leq x\}=\ell((0,x]\cap\tilde\Gamma),$ thus it defines a measure with  density function $I_{\tilde\Gamma}(t),$ and, by the Lebesgue differentiation theorem, 
\begin{equation}\label{Lebes}
 \lim_{h\to 0}\frac{F_G(t+h)-t}h=1 \ \mbox{ for almost every } t\in \tilde\Gamma.
\end{equation}
Now, from
\begin{eqnarray}\label{ellimite}
\lefteqn{\ell\big\{t\in \tilde\Gamma:{\textstyle t +\frac{B_{n}^F(t)}{\sqrt{n}} }>F_G({\textstyle t +\frac{B_{m}^G(t)}{\sqrt{m}} })+C_{n,m}\big\}}
\\ \nonumber
&=&\ell\big\{t\in \tilde\Gamma: {\textstyle \sqrt{\frac{m+n}n}}{\textstyle {B_{n}^F(t)} }>{\textstyle \sqrt{\frac{m+n}m}}B_{m}^G(t)(F_G({\textstyle t +\frac{B_{m}^G(t)}{\sqrt{m}} })-t)/{\textstyle \frac{B_{m}^G(t)}{\sqrt{m}}}+{\textstyle\sqrt{m+n}C_{n,m}}\big\}\\ \nonumber
 &\stackrel{d}{=}&\ell\big\{t\in \tilde\Gamma: {\textstyle \sqrt{\frac{m+n}n}}{\textstyle {B^F(t)} }>{\textstyle \sqrt{\frac{m+n}m}}B^G(t)(F_G({\textstyle t +\frac{B^G(t)}{\sqrt{m}} })-t)/{\textstyle \frac{B^G(t)}{\sqrt{m}}}+{\textstyle\sqrt{m+n}C_{n,m}}\big\},
\end{eqnarray}
where $B_F$ and $B_G$ are independent standard Brownian bridges, \eqref{Lebes}, the expression of $C_{n,m}$ and dominated convergence imply convergence to
$$
\ell\big\{t\in \tilde\Gamma, \lambda^{-1/2}{B_F(t)} -(1-\lambda)^{-1/2}{B_G(t)}>0\big\}.
$$ 
Finally, independence between $B_F$ and $B_G$ gives that $\lambda^{-1/2}{B_F(t)} -(1-\lambda)^{-1/2}{B_G(t)}$ is a scaled Brownian bridge (it can be written as $(\lambda^{-1}+(1-\lambda)^{-1})^{1/2}{B(t)}$, where $B(t)$ is a standard Brownian bridge). Therefore the limit law  in (\ref{ellimite})  is that  $\ell\{t \in \tilde\Gamma: B(t)>0\}.$
\FIN

\begin{Nota}{\rm
It is obvious that, for any Borel set $A\subset [0,1]$, and a Brownian bridge $B$, the distribution of $\ell( \{B>0\}\cap A)$ is supported by $[0,\ell(A)]$. One could conjucture that this distribution should be also uniform on $(0,\ell(A))$. However, a second thought shows that this distribution, in fact, depends on  the set $A$ and that {
it could even be non-continuous.  It is well known (see e.g. pag. 42 in \cite{Shorack}) that $P(B(t)\neq 0 \mbox{ for   }a<t<b)\neq 0$ if $0<a<b<1,$ thus     if  $A$ is contained in  $[a,b]$,} then the probability of the event $\left\{\ell( \{B>0\}\cap A) = \ell(A)\right\}$ is strictly positive. In fact, this distribution has two atoms:  at $\ell(A)$ and at 0.
\FIN}
\end{Nota}

\section{Rates of convergence}\label{puntos_separados}

When the set $\tilde{\Gamma}$ is negligible, Theorems \ref{one-side} and \ref{Theo_Noconsistency} yield convergence of Galton's rank statistic to the index $\gamma(F,G)$. We investigate in this section the rate of convergence in this result when the contact set $\Gamma^*$ (recall Definition \ref{Defin.ContactPoint}) is finite. The following simple result will be crucial in our analysis.
\begin{Lemm}\label{localizacion}
Assume that $[a,b]\subset [0,1]\backslash \Gamma^*$ is such that $t-F_G(t)>\delta>0$ for every $t\in[a,b]$. If  $\frac{n}{n+m}\to \lambda\in(0,1)$ then, for every $\varepsilon>0$ such that $a+\varepsilon<b-\varepsilon$, we have 
a.s. eventually
$$
[a+\varepsilon,b-\varepsilon]
=
\big\{t\in[a+\varepsilon,b-\varepsilon]:\, F_n^{-1}(t)>G_m^{-1}(t)\big\}
=
\big\{t\in[a+\varepsilon,b-\varepsilon]:\, F^{-1}(t)>G^{-1}(t)\}.$$ 
The same conclusion holds if $t-F_G(t)<-\delta$ for every $t\in[a,b]$.
\end{Lemm}

\noindent \textbf{Proof:} We have $F^{-1}(t)>G^{-1}(t)$ for every $t\in [a,b]$. 
Using the representation \eqref{Eq.represent}, 
$$
\big\{t\in[a+\varepsilon,b-\varepsilon]:\, F_n^{-1}(t)>G_m^{-1}(t)\big\}
=
\big\{t\in[a+\varepsilon,b-\varepsilon]:\, t+\frac{u_n(t)}{\sqrt{n}}>F_G(t+\frac{v_m(t)}{\sqrt{m}})\big\}.
$$
 Without loss of generality we can assume that the chosen version of $u_n$ satisfies $\sup_{0\leq t\leq 1}|u_n(t)|$ is a.s. {bounded}, and the same for $v_m$. Then, a.s., we have that for all $t\in [a+\varepsilon,b-\varepsilon]$, eventually $t+\frac{v_m(t)}{\sqrt{m}}\in[a,b]$ and therefore
$F_G\big(t+\frac{v_m(t)}{\sqrt{m}} \big)<t+\frac{v_m(t)}{\sqrt{m}}-\delta<t+\frac{u_n(t)}{\sqrt{n}}$ for large enough $n$ and $m$ and the result follows. 
The same argument fixes the case $t-F_G(t)<-\delta$.
\FIN

Now, recalling (see \eqref{Eq.Statistic.2}) the notation 
$$\ell_{n,m}^{t_0} :=
\ell\left( \{ F_n^{-1}>G_m^{-1}\} \cap(t_0-\eta,t_0+\eta)\right)
-
\ell\left( \{F^{-1}>G^{-1}\}\cap(t_0-\eta,t_0+\eta)\right),$$ we obtain, as an inmediate consequence of Lemma \ref{localizacion} and Lemma \ref{Lemm.SignoConstante} and the subsequent comments, the following result.

\begin{Coro}\label{localizacion2}
If $\Gamma^*=\{t_1,\ldots,t_k\}$ , $k>0$, $\frac{n}{n+m}\to \lambda\in(0,1)$ and $\eta>0$ is such that $\{t_i\}=\Gamma^*\cap (t_i-\eta,t_i+\eta)$, $i=1,\ldots,k$, then for $s>0$
$$n^{s}(\gamma(F_n,G_m)-\gamma(F,G))=n^s\sum_{i=1}^k \ell_{n,m}^{t_i}+o_P(1).$$
\end{Coro}

The main consequence of Lemma \ref{localizacion} and Corollary \ref{localizacion2} is that when $\Gamma^*$ is finite the key to the asymptotic behaviour of $\gamma(F_n,G_m)$ is the (joint) asymptotic behaviour of $\ell_{n,m}^{t_i}$. We address this problem in this section when $\Gamma^*$ consists of regular contact points. We note that these regular contact points (recall \eqref{Eq.DesarrolloDelta}) are elements of $\tilde{\Gamma}^*_F$. This, apparently, excludes contact points  in $\tilde{\Gamma}^*_G$ but not in $\tilde{\Gamma}^*_F$ or points which would be regular if we exchange the roles of $F$ and $G$ but 
are not with the present definition.
However, these cases can often be handled with the same approach. To see this, observe that when $\Gamma^*$ is finite (recall the concluding remarks in Section \ref{Sec.Notation}) we have that $\ell(t\in A:F^{-1}(t)\leq G^{-1}(t))=\ell(t\in A:F^{-1}(t)< G^{-1}(t))$ for every measurable $A$. 

If we assume further that $F$ and $G$ have no common discontinuity point (see Proposition \ref{common.disc} and the more general Proposition \ref{ultima}, involving just local conditions,} in the Appendix), then $\ell(t\in A:F_n^{-1}(t)\leq G_m^{-1}(t))=\ell(t\in A:F_n^{-1}(t)< G_m^{-1}(t))$ a.s. and we see that
\begin{eqnarray*}
\ell_{n,m}^{t_0} &=&-\big(
\ell\left( \{ F_n^{-1}<G_m^{-1}\} \cap(t_0-\eta,t_0+\eta)\right)
-
\ell\left( \{F^{-1}<G^{-1}\}\cap(t_0-\eta,t_0+\eta)\right)\big)\\
&=&:-\tilde{\ell}_{m,n}^{t_0}\quad \mbox{a.s.}.
\end{eqnarray*}
Observe that $\tilde{\ell}_{m,n}^{t_0}$ is the same statistic as ${\ell}_{n,m}^{t_0}$ after exchanging the roles of the $X$ and the $Y$ samples.
Hence, we restrict our analysis to points in $\Gamma^*_F$. Our results hold for points in $\tilde{\Gamma}^*_G$ with obvious changes.

We note that for every regular
contact point, $t_0$, there exists  $\eta^*>0$ such that
\begin{equation}\label{Eq.IS}
\mbox{$\mbox{sgn}(F_G(t)-t)$ is non-null and constant on each of $(t_0-\eta^*,t_0)$ and $(t_0,t_0+\eta^*)$.
}
\end{equation}

We recall from the final comments in Section \ref{Sec.Notation} that, by taking $\eta^*$ small enough (to exclude other contact points from the interval),  $\mbox{sgn}(F_G(t)-t)=\sign(G^{-1}(t)-F^{-1}(t))$ for every $t\in(t_0-\eta^*,t_0)\cup(t_0,t_0+\eta^*)$. Now, if \eqref{Eq.IS} holds, the study of 
$\ell_{n,m}^{t_0}$ can be carried out through the study, for $\eta \in (0,\eta^*)$, of the pieces 
\[
{L_{n,m}^>:=}
\int_{t_0-\eta}^{t_0}I_{\left\{F_n^{-1}(s) > G_m^{-1}(s) \right\}}ds \ \ 
\mbox{ and } 
\ \ 
{R_{n,m}^>:=} \int_{t_0}^{t_0+\eta}I_{\left\{F_n^{-1}(s) > G_m^{-1}(s) \right\}}ds, 
\]
\[
{L_{n,m}^<:=} \int_{t_0-\eta}^{t_0}I_{\left\{F_n^{-1}(s) \leq G_m^{-1}(s) \right\}}ds \ \ 
\mbox{ and } \ \ 
{R_{n,m}^<:=}\int_{t_0}^{t_0+\eta}I_{\left\{F_n^{-1}(s) \leq G_m^{-1}(s) \right\}}ds,
\]
 corresponding to the interval(s) (if any) where $F^{-1}>G^{-1}$.  For example, for a crossing point $t_0$ such that $F^{-1}<G^{-1}$ on $(t_0-\eta, t_0)$ and $F^{-1}>G^{-1}$ on $(t_0, t_0+\eta)$, $\ell_{n,m}^{t_0}=L_{n,m}^>-R_{n,m}^<$ (recall that for small enough $\eta$ this happens when $C_L(t_0)>0,C_R(t_0)<0$). With this notation we are ready for the proof of Theorem \ref{pieces}.

\medskip
\noindent \textbf{Proof of Theorem \ref{pieces}.} We assume, for instance, that $C_L>0$, and $r_L\geq r_R$, thus $r=r_L$. The other cases can be handled similarly.
We note that $\ell_{n,m}^{t_0}=L_{n,m}^>+R_{n,m}^>$ if $C_R>0$, while $\ell_{n,m}^{t_0}=L_{n,m}^>-R_{n,m}^<$ if $C_R<0$.
We consider first the case $r_L>1$. We set $d_n={\left({n+m}\right)}^{1/2r}$ and prove next that 
\begin{equation}\label{aux1}
d_n L_{n,m}^>\convw \ell
\left\{y< 0:B_\lambda({t_0}) > C_L\, |y|^r 
\right\}.
\end{equation}
To check this we note that, using \eqref{eventualbound2}, \eqref{Eq.Las Tildes}, 
Lemma \ref{Lem.BaseDeLoDeTasio} and \eqref{h-func}, it is enough to prove that
\begin{equation}\label{Eq.Cond.2}
d_n\ell\big\{t\in \mathcal I:{\textstyle t +\frac{B_1(t)}{\sqrt{n}} }>F_G({\textstyle t +\frac{B_2(t)}{\sqrt{m}} })-C_{n,m}\big\} \convw
\ell
\left\{y< 0:B_\lambda({t_0}) > C_L\, |y|^r 
\right\} 
\end{equation}
and  similarly with 
$
d_n\ell\big\{t\in \mathcal I:{\textstyle t +\frac{B_1(t)}{\sqrt{n}} }>F_G({\textstyle t +\frac{B_2(t)}{\sqrt{m}} })+C_{n,m}\big\} 
$, where $\mathcal I=[t_0-\eta,t_0]$.
The proofs are similar, hence,  we only prove  \eqref{Eq.Cond.2}. To ease notation we write $C(t)$ for $C_L$ when $t<0$ and for $C_R$ when $t>0$ and, similarly, $|t|^r$ will mean $|t|^{r_L}$ or $|t|^{r_R,}$ whenever $t<0$ or $t>0$. Then
\begin{eqnarray}
\nonumber
 I_{
\left\{t\in  \mathcal I: \;  t +\frac{B_1(t)}{\sqrt{n}} >
F_G\left( t +\frac{B_2(t)}{\sqrt{m}}\right)
-C_{n,m}  \right\}}
&=&
I_{
\left\{t\in  \mathcal I: \;  t +\frac{B_1(t)}{\sqrt{n}} >
t_0
+ \xi_m
+
C\left(\xi_m\right) \left|\xi_m\right|^r  + o\left(  \left|\xi_m\right|^r \right)-C_{n,m}  \right\}
} \ \ \
\\
\label{Eq.DesSerie}
&=&
I_{
\left\{t\in  \mathcal I: \;  B_n^{1,2}(t) >
\sqrt{n+m}\left( C\left(\xi_m\right)  \left|\xi_m\right|^r  +o\left(  \left|\xi_m\right|^r \right)-C_{n,m}\right) \right\}
}, \ \ \ \ \
\end{eqnarray}
where $\alpha_n=\sqrt{(n+m)/n}$,  $\beta_m=\sqrt{(n+m)/m}$, $B^{1,2}_n(t) = \alpha_n{B_1(t)}
-
\beta_m {B_2(t)}$, $\xi_m=\big(t+\frac{B_2(t)}{\sqrt{m}}  -t_0\big)$ and we have used that ${t_0}=F_G({t_0})$.
Denoting $\xi_m^*(y) = \frac y {d_n} + \frac {B^G({t_0}+ \frac y {d_n})}{\sqrt m}$, the change of variable $t={t_0} + \frac y {d_n}$,   and \eqref{Eq.DesSerie} lead to
\begin{eqnarray}
\label{Eq.DesarrolloEnSerie_2}
\lefteqn{d_n\ell\big\{t\in \mathcal I:{\textstyle t +\frac{B_1(t)}{\sqrt{n}} }>F_G({\textstyle t +\frac{B_2(t)}{\sqrt{m}} })-C_{n,m}\big\}}
\\ \nonumber
&=&
\int_{-d_n\eta}^{0} 
I_{\big\{B^{1,2}_n({t_0} + \frac y {d_n})
>
C\big(
\xi_m^*(y)
\big)
\big| \frac{(n+m)^{1/2r}}{(n+m)^{1/2r_R}}y + \frac{(n+m)^{1/2r}}{\sqrt m}{B_2({t_0}+ \frac y {d_n})}  \big|^r  + \sqrt{n+m} \; \left(o\left(  \left|\xi_m^*(y)\right|^r \right)-C_{n,m}\right)\big\}}  dy.
\end{eqnarray}

Since the Brownian bridges have continuous trajectories with probability one, they are  bounded and a.s.:
\begin{eqnarray*}
\sup_{y \in [-\eta d_n,0]}
\xi_m^*(y)  
&\leq&
\frac{\sup_{x\in [0,1] }|B_2(t)| }{\sqrt m} \to 0.
\\
\inf_{y \in [-\eta d_n,0]} \xi_m^*(y)
&\geq& 
- \eta - \frac{\sup_{x\in [0,1]}|B_2(t)| }{\sqrt m} \to - \eta.
\end{eqnarray*}
Thus, eventually, for every $y < 0$, $
\xi_m^*(y)
 \in [-\eta^*,\eta^*]$ and $C\big(
\xi_m^*(y)
\big)\geq \min (|C_L|,|C_R|)  >0$.  

This, the fact that $B_1$ and $B_2$ are a.s.  bounded, and also that $\frac{(n+m)^{1/2r}}{(n+m)^{1/2r_R}}$ is either equal to  one or, else, goes to infinity, yield that,  a.s.,  the order of $\sqrt{n+m}|\xi_m^*(y)|^r$ is $|y|^r$ or higher. Finally, the definition of $C_{n,m}$   allows us to conclude that there exists $M>0$ (depending on the particular trajectory of the Brownian bridges) such that
$$
I_{\left\{y \in [-\eta d_n,0]:\ B^{1,2}_n(t + \frac y {d_n})
>
C(
\xi_m^*(y)
)\left|y+{\frac{(n+m)^{1/2r}}{\sqrt m}}{B^2({t_0}+ \frac y {d_n})}  \right|^r  + \sqrt{n+m}\; \left(o(|\xi_m^*(y)|^r)-C_{n,m}\right)\right\}}
\leq
I_{\{ -M \leq y \leq 0\}}.
$$

Now, if we fix $y<0$ such that $B_\lambda({t_0}) \neq  C_L \, |y|^{r}$, then, a.s., 
\begin{equation*}
\label{Eq.r>1}
I_{\big\{B^{1,2}_n({t_0} + \frac y {d_n})
>C(\xi_m^*(y))
\left|y+{\frac{(n+m)^{1/2r}}{\sqrt m}}{B^2({t_0}+ \frac y {d_n})}  \right|^r +\sqrt{n+m}\; \left(o(|\xi_m^*(y)|^{r})-C_{n,m}\right) \big\}}
\to
I_{
\big\{B_\lambda({t_0}) > C_L\, |y|^{r}
\big\}.
}
\end{equation*}
From here, dominated convergence yields  \eqref{Eq.Cond.2}, hence, as noted above, \eqref{aux1}.  We note that the limit in \eqref{aux1} equals
\begin{equation}\label{aux2}
\sign(C_L)\Big({\textstyle \frac{(B_\lambda(t_0))^{{\mbox{\tiny sgn}(C_L)}}}{|C_L|}}\Big)^{1/r_0}. 
\end{equation}
A completely similar analysis shows that 
$$d_nR_{n,m}^>\convw  \sign(C_R)\Big({\textstyle \frac{(B_\lambda(t_0))^{{\mbox{\tiny sgn}(C_R)}}}{|C_R|}}\Big)^{1/r_0}I(r_R=r_0)$$
when $C_R>0$ ($d_nR_{n,m}^>$ vanishes in probability if $r_R<r_L=r_0$). Furthermore, we are using the same strong approximation to handle $d_n L_{n,m}^>$ and $d_nR_{n,m}^>$, which implies that there is weak convergence of $(d_n L_{n,m}^>, d_nR_{n,m}^>)$ and, consequently, of $d_n \ell_{n,m}^{t_0}=d_n(L_{n,m}^>+R_{n,m}^>)$. This completes the proof in the case $r_L\geq r_R$, $C_L>0,C_R>0$. The other cases with $r>1$ follow similarly.

The proof for the case $r_L=r_R=1$  goes along the same lines, the only difference being that, a.s., if $y<0$ satisfies  that $B_\lambda({t_0}) \neq  C_* \, |y+ \frac{B_2({t_0})}{ \sqrt{1-\lambda}}|$, where $*=L \mbox{ or } R$ whenever $\sign \left(y+ \frac{B_2({t_0})}{ \sqrt{1-\lambda}}\right)=-1 \mbox{ or } +1$,  we would have
$$
I_{\big\{B^{1,2}_n({t_0} + \frac y {d_n})
>C(\xi_m^*(y))
\left|y+{\frac{(n+m)^{1/2}}{\sqrt m}}{B_2({t_0}+ \frac y {d_n})}  \right|+\sqrt{n+m}\; \left(o(|\xi_m^*(y)|)-C_{n,m}\right) \big\}}
\to
I_{
\big\{B_\lambda({t_0}) > C_*\, |y+ \frac{B_2({t_0})}{ \sqrt{1-\lambda}}|
\big\}
}$$
and, by dominated convergence,
\begin{eqnarray*}
d_n L_{n,m}^> &\convw & \ell\big\{y<0:\, B_\lambda({t_0}) > C_*\, |y+ {\textstyle \frac{B_2({t_0})}{ \sqrt{1-\lambda}}}|
\big\}\\
&=& \ell\big\{y<0:\, B_\lambda({t_0}) >- C_L\, \big(y+ {\textstyle \frac{B_2({t_0})}{ \sqrt{1-\lambda}}}\big),y+ {\textstyle \frac{B_2({t_0})}{ \sqrt{1-\lambda}}}<0 \big\}\\
&&+ \ell\big\{y<0:\, B_\lambda({t_0}) > C_R\, \big(y+ {\textstyle \frac{B_2({t_0})}{ \sqrt{1-\lambda}}}\big),y+ {\textstyle \frac{B_2({t_0})}{ \sqrt{1-\lambda}}}>0 \big\}.
\end{eqnarray*}
The right side of the interval is dealt with in a similar way. In the case $C_L>0,C_R>0$, we end up concluding that
\begin{eqnarray*}
d_n R_{n,m}^> &\convw & \ell\big\{y>0:\, B_\lambda({t_0}) >- C_L\, \big(y+ {\textstyle \frac{B_2({t_0})}{ \sqrt{1-\lambda}}}\big),y+ {\textstyle \frac{B_2({t_0})}{ \sqrt{1-\lambda}}}<0 \big\}\\
&&+ \ell\big\{y>0:\, B_\lambda({t_0}) > C_R\, \big(y+ {\textstyle \frac{B_2({t_0})}{ \sqrt{1-\lambda}}}\big),y+ {\textstyle \frac{B_2({t_0})}{ \sqrt{1-\lambda}}}>0 \big\}.
\end{eqnarray*}
Hence,
\begin{eqnarray*}
d_n  \ell_{n,m}^{t_0}=d_n(L_{n,m}^>+R_{n,m}^>)&\convw & \ell\big(y: B_\lambda(t_0)>-C_L(y+{\textstyle \frac{B_2(t_0)}{\sqrt{1-\lambda}}}), y+{\textstyle \frac{B_2(t_0)}{\sqrt{1-\lambda}}}<0\big)\\
&&+
\ell\big(y: B_\lambda(t_0)>C_R(y+{\textstyle \frac{B_2(t_0)}{\sqrt{1-\lambda}}}), y+{\textstyle \frac{B_2(t_0)}{\sqrt{1-\lambda}}}>0\big)\\
&=&
{\textstyle \frac{(B_\lambda(t_0))^-}{C_L}}
+
{\textstyle \frac{(B_\lambda(t_0))+}{C_R}}=T_{1,1}(t_0;C_L,C_R).
\end{eqnarray*}
If $C_L>0,C_R<0$ we get
\begin{eqnarray*}
d_n R_{n,m}^< &\convw & \ell\big\{y>0:\, B_\lambda({t_0}) <- C_L\, \big(y+ {\textstyle \frac{B_2({t_0})}{ \sqrt{1-\lambda}}}\big),y+ {\textstyle \frac{B_2({t_0})}{ \sqrt{1-\lambda}}}<0 \big\}\\
&&+ \ell\big\{y>0:\, B_\lambda({t_0}) < C_R\, \big(y+ {\textstyle \frac{B_2({t_0})}{ \sqrt{1-\lambda}}}\big),y+ {\textstyle \frac{B_2({t_0})}{ \sqrt{1-\lambda}}}>0 \big\}.
\end{eqnarray*}
Therefore, 
\begin{eqnarray*}
d_n  \ell_{n,m}^{t_0}=d_n(L_{n,m}^>-R_{n,m}^>)&\convw & \ell\big(y: B_\lambda(t_0)>-C_L(y+{\textstyle \frac{B_2(t_0)}{\sqrt{1-\lambda}}}), y+{\textstyle \frac{B_2(t_0)}{\sqrt{1-\lambda}}}<0\big)\\
&&-\ell\big(y>0:\, y+{\textstyle \frac{B_2(t_0)}{\sqrt{1-\lambda}}}<0\big)\\
&&-\ell\big(y: B_\lambda(t_0)<C_R(y+{\textstyle \frac{B_2(t_0)}{\sqrt{1-\lambda}}}), y+{\textstyle \frac{B_2(t_0)}{\sqrt{1-\lambda}}}>0\big)\\
&&+\ell\big(y<0:\, y+{\textstyle \frac{B_2(t_0)}{\sqrt{1-\lambda}}}>0\big)\\
&=&{\textstyle \frac{(B_\lambda(t_0))^+}{C_L}-\frac{(B_2(t_0))^-}{\sqrt{1-\lambda}} +\frac{(B_\lambda(t_0))-}{C_R}+\frac{(B_2(t_0))^+}
{\sqrt{1-\lambda}}}\\
&=&{\textstyle \frac{(B_\lambda(t_0))^+}{C_L}+\frac{(B_\lambda(t_0))^-}{C_R}+\frac{B_2(t_0)}
{\sqrt{1-\lambda}}}=T_{1,1,}(t_0;C_L,C_R).
\end{eqnarray*}
The remaining cases are completely similar. We omit further details.
\FIN

\bigskip
Some comments are in order here. 
First note that, by focusing on the transform $F_G$, Theorem \ref{pieces} is able to handle virtual contact points for $F^{-1}$ and $G^{-1}$. As an illustration of this claim, assume $F^{-1}(t_0)<G^{-1}(t_0)\leq G^{-1}(t_0+)< F^{-1}(t_0+)$ ($t_0$ is then a virtual crossing point). 
As noted above, $F_G(t)=t_0$ in an interval $(t_0-\eta,t_0+\eta)$ for $\eta$ small enough, and \eqref{Eq.Delta} holds with  $r_R=r_L=1$, $C_R(t_0)=-1$ and $C_L(t_0)=+1$. Thus,  Theorem \ref{pieces} applies and gives (\ref{virtual.crossing}) below.

The case $F^{-1}(t_0)<G^{-1}(t_0)\leq F^{-1}(t_0+)< G^{-1}(t_0+)$ (a virtual tangency point) can be handled similarly, although it does not fit exactly in the setup of Theorem \ref{pieces}. In this case we have that, for some small enough $\eta, \delta>0$, $F_G(t)=t_0$, $t\in(t_0-\eta,t_0)$, $F_G(t)>t+\delta$, $t\in(t_0,t_0+\eta)$.
It is easy to see that, eventually,
$\ell_{n,m}^{t_0}=\ell\{t\in (t_0-\eta,t_0+\eta):\, t+\frac{u_n(t)}{\sqrt{n}}> t_0, t+\frac{v_m(t)}{\sqrt{m}}< t_0 \}$. From this point one can argue as in the proof of Theorem \ref{pieces} to obtain (\ref{tangencia.virtual}) below.

We include in the following proposition these results for virtual contact points. Notice that this proposition includes the possibility of non-continuous d.f.'s $F$ or $G$.

\begin{Prop}\label{Theo_Sing_disc}
Let  $t_0 \in \Gamma^* \cap (0,1)$, such that for some $\eta_0>0$, $(t_0-\eta_0,t_0+\eta_0)\cap \Gamma^*=\{t_0\}$. Then, for every small enough $\eta>0$, if $n/(n+m) \to \lambda \in (0,1)$ as $n,m\to \infty$, we  have that: 
\begin{enumerate}
\item[(i)]
{\em({\bf virtual crossing points})}
 If $F^{-1}(t_0)<G^{-1}(t_0)\leq G^{-1}(t_0+)< F^{-1}(t_0+)$, then, 
 \begin{equation}\label{virtual.crossing}
\sqrt{n+m} \ell_{n,m}^{t_0}\convw {\textstyle\frac{B_1(t_0)}{\sqrt{\lambda}}}.
\end{equation}

\item[(ii)] 
{\em({\bf  virtual tangency points})} If 
$F^{-1}(t_0)<G^{-1}(t_0)\leq F^{-1}(t_0+)< G^{-1}(t_0+)$,
then, 
\begin{equation}\label{tangencia.virtual}
\sqrt{n+m}\ell_{n,m}^{t_0}\convw {\textstyle \ell\{y:\, -\frac{B_2(t_0)}{\sqrt{1-\lambda}}>y,-\frac{B_1(t_0)}{\sqrt{\lambda}}<y  \}=\Big(\frac{B_1(t_0)}{\sqrt{\lambda}}- \frac{B_2(t_0)}{\sqrt{1-\lambda}}\Big)^+}.
\end{equation}
\end{enumerate}
\end{Prop}

We can  easily adapt Proposition \ref{Theo_Sing_disc} to the case $G^{-1}(t_0)<F^{-1}(t_0)\leq F^{-1}(t_0+)< G^{-1}(t_0+)$. In this case $F_G(t_0)<t_0<F_G(t_0+)$ (but $G_F(t)=t_0$ for $t$ close to $t_0$). We call this kind of virtual crossing a \textit{vertical} crossing, while we will refer to case \textit{(i)} as a \textit{horizontal} crossing. For vertical crossing points the argument above yields $\sqrt{n+m} \ell_{n,m}^{t_0}\convw -{\textstyle\frac{B_2(t_0)}{\sqrt{1-\lambda}}}$. Also, \textit{(ii)} corresponds to an \textit{upper} tangency point, in the sense that $F_G$ touches the identity at $t_0$ but remains above it in $(t_0-\eta,t_0+\eta)\backslash \{t_0\}$. With obvious changes we can deal with \textit{lower} tangency points, obtaining then 
$\sqrt{n+m}\ell_{n,m}^{t_0}\convw {\textstyle -\Big(\frac{B_1(t_0)}{\sqrt{\lambda}}- \frac{B_2(t_0)}{\sqrt{1-\lambda}}\Big)^-}$.

A further observation is that, since any contact order $r_0\geq 1$ is possible (see Example \ref{ejemplo.alternativo2}), we can obtain any rate of convergence $(m+n)^{-s}$, $s\leq1/2$ for $\ell^{t_0}_{n,m}$. As previously mentioned, only the case $r_L(t_0)=r_R(t_0)=1$ and $C_L(t_0)=-C_R(t_0)$  leads to asymptotic normality.

Finally, we note that the limiting expressions become simpler under regularity. In fact, if $h(t)=F_G(t)-t$ is $r$ times differentiable with continuity at a point $t_0\in (0,1)$, such that $h(t_0)=0$ and with derivatives $h^{k)}(t_0)=0$, $k=1,\ldots,r-1$ and $h^{r)}(t_0)\ne 0$, then $t_0$ is an isolated contact point in the sense of \eqref{Eq.DesarrolloDelta} with $r_L(t_0)=r_R(t_0)=r$. For odd $r\geq 3$ we have $C_R(t_0)=-C_L(t_0)=\frac{h^{r)}(t_0)}{r!}$ and the conclusion in Theorem \ref{pieces} reads
\begin{equation}\label{k.odd}
(n+m)^{\frac 1 {2r}}\ell_{n,m}^{t_0}\convw {\textstyle \Big(\frac{r!}{|h^{r)}(t_0)|}\Big)}^{1/r}\big(((B_\lambda(t_0))^{1/r})^+-((B_\lambda(t_0))^{1/r})^-\big),
\end{equation}
while for $r=1$ it becomes
\begin{equation}\label{1.odd}
(n+m)^{\frac 1 {2}}\ell_{n,m}^{t_0}\convw \sign(h'(t_0))\big({\textstyle \frac 1 {\sqrt{\lambda }}\frac{1}{h'(t_0)}B_1(t_0)+\frac 1 {\sqrt{1-\lambda }}\big(1+\frac{1}{h'(0)}\big)B_2(t_0)}\big). 
\end{equation}
For even $r\geq 2$ we have $C_R(t_0)=C_L(t_0)=\frac{h^{r)}(t_0)}{r!}$ and Theorem \ref{pieces} yields
\begin{equation}\label{k.even}
(n+m)^{\frac 1 {2r}}\ell_{n,m}^{t_0}\convw \sign(h^{r)}(t_0))2{\textstyle \Big(\frac{r!}{|h^{r)}(t_0)|}}\Big)^{1/r}\big((B_\lambda(t_0))^{\mbox{ \scriptsize sgn}( h^{r)}(t_0))}\big)^{1/r}.
\end{equation}

When the contact point is extremal, that is, when $t_0\in\{0,1\}$, the limiting r.v.'s in Theorem \ref{pieces} vanish. 
We prove now Theorem \ref{Theo_0_1.1}, showing that in this case there is weak convergence, at a faster rate, to a nondegenerate limiting distribution.  

\medskip
\noindent \textbf{Proof of Theorem \ref{Theo_0_1.1}.} 
Let us take $t_0=0$ and $C>0$. The cases with $C<0$ and/or $t_0=1$ are similar.
We handle first the case $r=1$. Then for small enough $\eta$ we have
$\ell\left( \{F^{-1}>G^{-1}\}\cap(0,\eta)\right)=0$. We recall that
\begin{equation} \label{Eq.SingularAt0_1}
\ell_{n,m}^0=  \ell \left( \{ F_n^{-1} > G_m^{-1}\} \cap(0,\eta)\right).
\end{equation}
We use the well-known fact that the joint law of
$
\frac 1 {S^1_{n+1}}(S^1_{1},\ldots,S^1_{n})
$
is the same as that the ordered sample of size $n$ of i.i.d. $U(0,1)$ r.v.'s. 
Thus,
\begin{equation} \label{Eq.Cruce.Cero_1}
(X_{(1)},\ldots,X_{(n)})\overset d = \Big( F^{-1}\big({\textstyle \frac{S^1_{1}} {S^1_{n+1}}}\big),\ldots,F^{-1}({\textstyle \frac{S^1_{n}} {S^1_{n+1}}}) \Big),
\end{equation}
with a similar expression for the $Y$-sample. 
From  (\ref{Eq.SingularAt0_1}) and (\ref{Eq.Cruce.Cero_1}) we see that
\begin{eqnarray*}
\nonumber
\ell_{n,m}^0
&\overset d = &
\int_0^\eta I_{\left\{
F^{-1}\big(\frac{S^1_{\lceil n t \rceil}} {S^1_{n+1}}\big)
>
G^{-1}\big(\frac{S^2_{\lceil m t \rceil}}{S^2_{m+1}}\big)  \right\}} dt
=
\frac 1 {n+m}\int_0^{(n+m)\eta}
I_{\big\{ 
F^{-1}(\xi_n^1(y))
>
G^{-1}(\xi_m^2(y))\big\} } dy,
\end{eqnarray*}
where $\xi_n^1(y) := S^1_{\lceil { \frac{n}{n+m}} y \rceil} /S^1_{n+1}$, and $\xi_m^2(y):=S^2_{\lceil { \frac{m}{n+m}} y \rceil} /S^2_{m+1}$.
Now \eqref{Eq.DesarrolloDelta} yields that,  if $y \in (0,(n+m)\eta)$,   then
\begin{eqnarray}
\nonumber
I_{\left\{
F^{-1}\left(\xi_n^1(y)\right)
>
G^{-1}\left(\xi_m^2(y)\right)
\right\} }
&=&
I_{\left\{
\xi_n^1(y)
>
F_G\left(\xi_m^2(y)\right)
\right\}}
\\
\label{Eq.Cruce.Cero_7}
&=&
I_{\left\{
{ {  (n+m)}}\xi_n^1(y)
>
(1+C){  (n+m)} {\xi_m^2(y)} + {  (n+m)}o({\xi_m^2(y)} )
\right\}}.
\end{eqnarray}

The SLLN  implies that there exists $\Omega_0$, with $P(\Omega_0)=1$, such that for every $\omega\in \Omega_0$, 
$
\frac{S^1_{n}}{n} \conv 1
\mbox{ and }
\frac{S^2_{m}}{m} \conv 1.
$ 
Therefore, for any $\delta^*>0$, if $\omega \in \Omega_0$, eventually
\begin{equation} \label{Eq.Cruce.0_4}
0 < \xi_m^2(y) \leq \frac{S^2_{\lceil { (1-\lambda +\delta^*)y \rceil} +1}}{S^2_{m+1}} \to 0.
\end{equation}
This and (\ref{Eq.Cruce.Cero_7}) show that if $\omega \in \Omega_0$,
\[
I_{\big\{
F^{-1}\big(\xi_n^1(y)\big)
>
G^{-1}\big(\xi_m^2(y)\big)
\big\} }
\conv
I_{\big\{
{ (1-\lambda)} S^1_{\lceil \lambda   y \rceil}
>
  \lambda  {
(1+C)}S^2_{\lceil (1-\lambda)  y \rceil}
\big\}},
\]
for every $y$ not belonging to the countable set $\{\frac j {1-\lambda}: \, j=0,1,\ldots\}\cup \{\frac j {\lambda}: \, j=0,1,\ldots\}$. Clearly, in $\Omega_0$ we have
$$\lim_{y\to\infty} \frac{(1-\lambda) S^1_{\lceil \lambda y \rceil}}{\lambda S^2_{\lceil (1-\lambda) y \rceil}}=1.$$
Hence, the fact that $C>0$ gives that, in $\Omega_0$, 
$I_{\left\{   (1-\lambda)S^1_{\lceil \lambda   y \rceil}
>
{
\lambda(1+C)}  S^2_{\lceil (1-\lambda)  y \rceil}\right\} }=0$ for large enough $y$. This shows that
$
\int_0^\infty I_{\left\{
{ (1-\lambda)} S^1_{\lceil \lambda   y \rceil}
>
{  \lambda } {
(1+C)} S^2_{\lceil (1-\lambda)  y \rceil}
\right\}}
dy
$ 
is an a.s. finite r.v.. 

We will conclude (\ref{distr_limite}) as soon as we prove that for every  $\omega \in \Omega_0$ we can apply dominated convergence. To check this, notice that (\ref{Eq.Cruce.0_4}) gives that, for $m$ large enough, 
\begin{eqnarray}
 \label{Eq.Cruce.Cero_4}
\lefteqn{I_{\left\{
 (m+n) \xi_n^1(y)
>
   (m+n) {
(1+C)} \xi_m^2(y) + (m+n) o(\xi_m^2(y))
\right\}
}
}
\\
\nonumber
 &\hspace*{25mm}\leq &
I_{\left\{
(m+n)\xi_n^1(y)
>
(1+C/2)  (m+n)  \xi_m^2(y) 
\right\}
}.
\end{eqnarray}

Now, for every $\omega \in \Omega_0$,  there exist a natural number and a positive real number depending on $\omega$, $N(\omega)$ and $Y(\omega)$, such that, if $n \ge N(\omega)$   then both $S^1_{n+1}/n$ and $S^2_{m+1}/m$ are close to one, and, if we take $y \geq Y(\omega)$, then, both $\frac{S^1_{\lceil { \frac{n}{n+m}} y \rceil}}{\frac n{n+m}}$ and  $\frac{S^2_{\lceil { \frac{m}{n+m}} y \rceil}}{\frac m{n+m}}$ are close to $y$. This completes the proof for the case $r=1$,  since (\ref{Eq.Cruce.Cero_4}) gives that, for all $n\geq N(\omega)$, 
 \[
I_{\left\{
\xi_n^1(y)
>
{
(1+C(\xi_m^2(y)))}   (m+n) \xi_m^2(y) +  (m+n) o(\xi_m^2(y))
\right\}
}
\leq I_{[0,Y(\omega)]}.
 \]

For the case $r>1$ we assume, again, $C>0$ and
observe that the r.v. $T_{r,r}(0;C,C)$
is a.s. finite  (this follows, for instance, from the fact that, a.s., $W_0(y)/y\to 0$ as $y\to \infty$). Now $\ell^0_{n,m}$ has the same expression as in (\ref{Eq.SingularAt0_1}). We will use the same notation as in Lemma \ref{Lem.BaseDeLoDeTasio}. First, we have that
\begin{eqnarray*}
\ell_{n,m}^0&=&
\ell\Big\{t \in (0,\eta):F^{-1}\big(
{\textstyle t + \frac{u_n(t)}{\sqrt n}}
\big) > G^{-1}\big({\textstyle 
t+ \frac{v_m(t)}{\sqrt m}}
\big) \Big\}
\\
&=&
\ell\Big\{t \in (0,\eta):{\textstyle t+ \frac{u_n(t)}{\sqrt n}}
 > F_G\big({\textstyle 
t + \frac{v_m(t)}{\sqrt m}}
\big) \Big\}.
\end{eqnarray*}
Therefore, if we take $f$ equal to the identity and $g=F_G$ in Lemma \ref{Lem.BaseDeLoDeTasio}, we only need to show that
\begin{eqnarray} \nonumber
\lefteqn{d_n 
\ell\left(\left\{ {\textstyle \frac{t+ B_n^F(t)}{\sqrt n}} >\tilde{F_G} \left({\textstyle t }\right)-L_n\right\}\cap(0,\eta)\right)}\hspace*{3cm}\\
\label{Eq.Condicion1.Theor4.6}
&\to_w &
\ell\big\{y\in (0,\infty):\, W_0(y)>(\lambda (1-\lambda))^{1/2}C(0)y^{r}\big\},
\end{eqnarray}
 and similarly for $d_n 
\ell\left(\left\{ {\textstyle t + \frac{B_n^F(t)}{\sqrt n}} >\tilde{F_G} \left({\textstyle t}\right)+L_n\right\}\cap(0,\eta)\right)$, where,
now, $d_n=(n+m)^{\frac 1 {2r-1}}$ (since $t+ \frac{B_m^G(t)}{\sqrt m}$ can take negative values, we take $F_G(t)= F_G(0)$ for $t<0$; notice that $t+ \frac{B_m^G(t)}{\sqrt m} \to t >0$, hence, this assumption has no effect in the limit) . The proofs are the same, thus we only consider (\ref{Eq.Condicion1.Theor4.6}).

We can assume, without loss of generality that $B_{n}^F(t)=d_n^{-1/2}(W_F(d_n t)-tW_F(d_n))$, and
$B_{m}^G(t)=d_n^{-1/2}(W_G(d_n t)-tW_G(d_n))$, $0\leq t\leq 1$ with $W_F,W_G$ independent Brownian motions. Thus, the change of variable $t=y/d_n$ and the fact that $\sqrt{(n+m)d_n}=d_n^{r}$ give 
\begin{eqnarray*}
\lefteqn{d_n\ell_{n,m}^0=
d_n \int_0^\eta I_{\Big(\alpha_nB_{n}^F(t)
>\beta_mB_{m}^G(t)+
\sqrt{n+m}{
C\big(t+{\textstyle \frac {B_{m}^G(t)}{\sqrt{m}}} \big)}\big|t+{\textstyle \frac {B_{m}^G(t)}{\sqrt{m}}} \big|^{r} +\sqrt{n+m}\; \big( o\big(\big|t+{\textstyle \frac {B_{m}^G(t)}{\sqrt{m}}} \big|^{r_{R}}\big)-L_{n,m}\big)
 \Big)}dt 
 }
\\
&=& 
\int_0^{d_n\eta} 
I_{\Big(\alpha_n\big(W_{F}(y)-y{\textstyle \frac{W_F(d_n)}{d_n}}\big)
>
\beta_m\big(W_{G}(y)-y{\textstyle \frac{W_G(d_n)}{d_n}}\big)
+ 
{
C((\xi_n(y))}d_n^{r} |\xi_n(y)|^{r}   +d_n^{r} \big(o(|\xi_n(y)|^{r})-L_{n,m}\big)\Big)
}dy,
\end{eqnarray*}
where $\alpha_n=((m+n)/n)^{1/2}$, $\beta_m=((m+n)/m)^{1/2}$ and $\xi_n(y)=\frac y {d_n} +\frac 1 {\sqrt m} B_{m}^G(\frac y {d_n})$.  

As it is well known, there exists $\Omega_0\in \sigma$, with $P(\Omega_0)=1$ such that, if $\omega \in \Omega_0$, then,  $W_i$ is continuous, $W_i(x)/x \to 0$, as $x \to \infty, i=F,G$ and  the set 
$$
\left\{y:\, \lambda^{-1/2}W_F(y) =(1-\lambda)^{-1/2}W_G(y)+C(0)y^{r}\right\}
$$
 has Lebesgue measure zero. If we fix $\omega \in \Omega_0$, then,  we have that
\[
\sup_{y \in [0,d_n\eta]} |\xi_n(y)| \le \eta +\frac 1{\sqrt{m d_n}} \sup_{y \in [0, d_n\eta]} \left|
W_F(y) - y \frac {W_F(d_n)}{d_n}
\right|
\to \eta,
\]
and we can conclude that, eventually, $\{\xi_n(y): y \in [0,d_n \eta]\} \subset [0,\eta^*]$, and, consequently, from an index onward,
$
\inf_{y \in [0,d_n\eta]}  C(\xi_n(y)) \geq \inf_{h \in [0,\eta^*]} |C(h)| >0.
$ 
 On the other hand, we have $d_n^{r} L_{n ,m}\to 0$ and
 \[
d_n^{r} |\xi_n(y)|^{r}
=
\left|
y + \beta_m{\textstyle \frac {W_{G}(y)-y{\textstyle \frac{W_G(d_n)}{d_n}}}{d_n^{r-1}}}
\right|^{r}
=y(1+o(1)) \to \infty \mbox{, as } y \to \infty.
\]
Therefore, there exists a constant $M$ (which possibly depends on the chosen $\omega$) such that 
$$
I_{\left(\alpha_n\big(W_{F}(y)-y{\textstyle \frac{W_F(d_n)}{d_n}}\big)
>
\beta_m\big(W_{G}(y)-y{\textstyle \frac{W_G(d_n)}{d_n}}\big)
+ 
C(\xi_n(y))d_n^{r} |\xi_n(y)|^{r} + d_n^{r} \big(o(|\xi_n(y)|^{r})-L_{n,m}\big)\right)}
\leq
 I_{\big\{0\leq y\leq M\big\}},
 $$
for every  large enough $n$. Moreover,
\begin{eqnarray*}
\lefteqn{I_{\left\{\alpha_n\left(W_{F}(y)-y{\textstyle \frac{W_F(d_n)}{d_n}}\right)
>
\beta_m\left(W_{G}(y)-y{\textstyle \frac{W_G(d_n)}{d_n}}\right)
+ 
C(\xi_n(y))d_n^{r} |\xi_n(y)|^{r}  + d_n^{r} \big(o(|\xi_n(y)|^{r})-L_{n,m}\big)\right\}}
}\hspace*{5cm}
\\
&\to& 
I_{\left\{\lambda^{-1/2}W_F(y)-(1-\lambda)^{-1/2}W_G(y)> C(0)y^{r}\right\}}.
\end{eqnarray*}

These observations allow to apply dominated convergence to conclude that, for this $\omega$,
\begin{eqnarray*}
\lefteqn{d_n 
\ell\big\{(\tilde{F}_n^{-1}>\tilde{G}_m^{-1}-L_{n,m})\cap(0,\eta)\big\}}
\hspace*{2cm}
\\
&\to&
\ell\big\{y\in (0,\infty):\, \lambda^{-1/2} W_F(y)-(1-\lambda)^{-1/2}W_G(y)
>
C(0)y^{r}\big\}
\end{eqnarray*}

The fact that $(\lambda (1-\lambda))^{1/2}(\lambda^{-1/2} W_F(y)-(1-\lambda)^{-1/2}W_G(y))$ is a standard Brownian motion yields \eqref{Eq.Condicion1.Theor4.6}. 
\FIN

We prove next Theorem \ref{TCL.principal}, a global asymptotic result for Galton's statistic under the assumption of a finite contact set consisting of regular contact points. From a technical point of view the main issue here is to prove asymptotic independence between the localized statistics around central and extremal contact points. 

\medskip
\noindent \textbf{Proof of Theorem \ref{TCL.principal}:} 
From Corollary \ref{localizacion2} it is enough to prove that $(n+m)^{\frac 1 {2r_0}}(\ell_{n,m}^{t_i})_{1\leq i\leq k}$ converges weakly. This follows trivially if $\Gamma^*\subset (0,1)$ after checking that the strong approximation used in the proof of Theorem \ref{pieces} allows to deal with all the $\ell_{n,m}^{t_i}$ simultaneously.
Hence, it suffices to prove asymptotic independence among $\ell_{n,m}^0$, $(\ell_{n,m}^{t_i})_{i:t_i\in(0,1)}$ and  $\ell_{n,m}^1$ when $0$ or $1$ (or both) are contact points. Let us assume, for instance, that $\Gamma^*=\{0<t_1\cdots<t_s<1 \}$ and set $A_n=(n+m)^{\frac 1 {2r_0}}\ell_{n,m}^0$, $B_n=(n+m)^{\frac 1 {2r_0}}(\ell_{n,m}^{t_i})_{1\leq i\leq s}$ and $C_n=(n+m)^{\frac 1 {2r_0}}\ell_{n,m}^1$. 
We have that there exist $A,B,C$ such that $A_n\convw A$, $B_n\convw B$, $C_n\convw C$. Assume $(\tilde{A},\tilde{B},\tilde{C})$ is a random vector with $\tilde{A},\tilde{B},\tilde{C}$ independent, $\tilde{A}\overset  d= A,\tilde{B}\overset  d= B$ and $\tilde{C}\overset  d= C$ and consider $(\tilde{A}_n,\tilde{B}_n,\tilde{C}_n)$, with the same properties with respect $(A_n,B_n,C_n)$. $A_n$ is a function of the smallest $\lceil \eta n \rceil$ elements in the $X$ sample and the smallest 
$\lceil \eta m \rceil$ elements in the $Y$ sample. Similarly, $B_n$ and $C_n$ are functions of the central and upper order statistis. If $d_{TV}$ denotes the distance in total variation, then there exists a universal constant $H>0$ such that
$$d_{TV}(\mathcal{L}(A_n,B_n,C_n),\mathcal{L}(\tilde{A}_n,\tilde{B}_n,\tilde{C}_n))\leq H \Big[\frac{\eta (1-t_1-\eta)}{t_1-2\eta}+ 
\frac{\eta (t_s+\eta)}{1-t_s-2\eta}\Big]^{1/2}$$
for small enough $\eta$ (this follows from Theorem 4.2.9 and Lemma 3.3.7 in \cite{Reiss1989}). If $\rho$ denotes the Prokhorov metric, then the fact that  $\rho(\mu_1,\mu_2)\leq d_{TV}(\mu_1,\mu_2)$ implies
\begin{equation}\nonumber
\rho(\mathcal{L}(A_n,B_n,C_n),\mathcal{L}(\tilde{A}_n,\tilde{B}_n,\tilde{C}_n))\leq H \Big[\frac{\eta (1-t_1-\eta)}{t_1-2\eta}+ 
\frac{\eta (t_s+\eta)}{1-t_s-2\eta}\Big]^{1/2}.
\end{equation}

We prove now that $(A_n,B_n,C_n)\convw (\tilde{A},\tilde{B},\tilde{C})$. Obviously $(\tilde{A}_n,\tilde{B}_n,\tilde{C}_n)\convw (\tilde{A},\tilde{B},\tilde{C})$. Having weakly convergent components, $(A_n,B_n,C_n)$ is tight. To complete the proof it suffices to show that for any weakly convergent subsequence $(A_{n'},B_{n'},C_{n'})\convw \gamma$, necessarily $\gamma=\mathcal{L} (\tilde{A},\tilde{B},\tilde{C})$. To check this, we observe that, since $\rho$ metrizes the weak convergence, we have 
\begin{equation}\label{asint.indep}
\rho(\gamma,\mathcal{L}(\tilde{A},\tilde{B},\tilde{C}))\leq H \Big[\frac{\eta (1-t_1-\eta)}{t_1-2\eta}+ 
\frac{\eta (t_s+\eta)}{1-t_s-2\eta}\Big]^{1/2}.
\end{equation}

Now, using Corollary \ref{localizacion2} we see that we can repeat the argument leading to \eqref{asint.indep} for every small enough $\eta$. Hence,
$\rho(\gamma,\mathcal{L}(\tilde{A},\tilde{B},\tilde{C}))=0$. This completes the proof.
\FIN

\subsection{Some examples and extensions}\label{examplessec}

We provide here some simple examples that illustrate the different limiting distributions for $\ell_{n,m}^{t_0}$ that result from
Theorems \ref{pieces} and \ref{Theo_0_1.1}. Later we give simple sufficient conditions under which extremes have no influence on the asymptotic behaviour of $\gamma(F_n,G_m)$ and give a simplified version of Theorem \ref{TCL.principal} under the assumption that $F$ and $G$ have regular densities (Theorem \ref{dosmuestras_2}). Finally, we consider the case of finitely supported distributions (Theorem \ref{finite.support}).

\begin{Ejem}\label{ejemplo.alternativo}
{\rm In this example $G(t)=t$ (the uniform law on $(0,1)$). For $r>0$ we consider the quantile function
$F^{-1}(t)=\frac 1 2 +\sign(t-\frac 1 2)|t-1/2|^r$, $0\leq t\leq 1$. Now we have 
$F(x)=\frac 1 2+\sign(x-\frac 1 2)|x-\frac 1 2|^{1/r}$, $\frac 1 2-\frac 1 {2^r}\leq x\leq \frac 1 2+\frac 1 {2^r}$,
$F_G=F$ and $F_G(\frac 1 2)=\frac 1 2$. Thus, $\frac 1 2$ is a contact point. If $r< 1$ then $F_G'(t)=\frac 1 r 
|t-\frac 1 2|^{\frac 1 r -1}$. 
{In particular, $F_G$ is Lipsichitz in a neighbourhood of $\frac 1 2$.} We easily check that
$\Delta(h)=-h+\sign(h)|h|^{1/r}=-h+o(h)$, that is, $\frac 1 2$ is an isolated {regular} contact point (a crossing point) with intensities
$r_L=r_R=1$ and constants $C_R=-C_L=-1$. We can apply Theorem \ref{pieces} to  conclude that 
$$(n+m)^{1/2}\ell_{n,m}^{t_0}\convw {\textstyle \frac{B_1(\frac 1 2)}{\sqrt{\lambda}}}.$$

If $r>1$ then $F_G'(\frac  1 2)=+\infty$ and { $F_G$ is not Lipschitz around} the contact point. However, following the reasoning after Corollary \ref{localizacion2}, we have that $\ell_{n,m}^{t_0}=-\tilde{\ell}_{m,n}^{t_0}$ and 
 we can handle this case exchanging the roles of the $F$ and $G$ samples and studying $G_F(t)=F^{-1}(t)$.
Now $G_F'(t)=r |t-\frac 1 2|^{r-1}$ and $G_F$ is Lipschitz in  a neighbourhood of $\frac 1 2$. Furthermore,
$\Delta(h)=-h+\sign(h)|h|^{r}=-h+o(h)$. Thus we can, again, apply Theorem \ref{pieces} to $\tilde{\ell}_{m,n}^{t_0}$ (with $r_L=r_R=1, C_R=-C_L=-1$) and 
conclude  that $(n+m)^{1/2}\tilde{\ell}_{m,n}^{t_0}\convw \frac{B_2(\frac 1 2)}{\sqrt{1-\lambda}}$. Hence, for $r>1$ we see that
$$(n+m)^{1/2}\ell_{n,m}^{t_0}\convw {\textstyle -\frac{B_2(\frac 1 2)}{\sqrt{1-\lambda}}}.$$
\FIN
}
\end{Ejem}

\begin{Ejem}\label{ejemplo.alternativo2}
{\rm Now $F$ denotes the d.f. of the  uniform law on $(0,1)$ and $G^{-1}(t)=t+\sign(t-\frac 1 2)|t-1/2|^r$, $0\leq t\leq 1$.
As before, $\frac 1 2$ is a contact point. For $r\geq 1$ $F_G=G^{-1}$ is differentiable, with $F_G'(t)=1+r |t-1/2|^{r-1}$. 
We have $\Delta(h)=\sign(h)|h|^r$, that is, Theorem \ref{pieces} can be applied here with $r_L=r_R=r$, $C_R=-C_L=1$. Thus, for $r=1$ we get 
$$(n+m)^{1/2}\ell_{n,m}^{t_0}\convw {\textstyle \frac{B_1(\frac 1 2)}{\sqrt{\lambda}}+ \frac{2B_2(\frac 1 2)}{\sqrt{1-\lambda}}},$$
while for $r>1$ we obtain 
$$(n+m)^{1/2r}\ell_{n,m}^{t_0}\convw {\textstyle ((B_\lambda(\frac 1 2))^+)^{1/r}- ((B_\lambda(\frac 1 2))^-)^{1/r}}.$$
The case $0<r<1$ can be handled exchanging the roles of the two samples, as in Example \ref{ejemplo.alternativo}. We omit details.
\FIN
}
\end{Ejem}

\begin{Ejem}\label{Ejem.t_nu}
{\rm
Here we consider a Student's $t$ location  model. Let $F=F_\nu$ be a $t$-distribution with $\nu>0$ degrees of freedom and $G(x)={G_\nu}(x)={F_\nu}(x-\mu)$, for some $\mu>0$. Obviously, in this case ${F}^{-1}(0)={G}^{-1}(0)=-\infty$ and $F_G(0)=0$.
We write $f_v$ for the density of $F_\nu$. To ease notation,  we set $s=(\nu +1)/2$, write $K$  for a non-null  generic constant which can change from line to line (in particular,  $f_\nu(t) = K(\nu + t^2)^{-s}$) and $f(x)\approx g(x)$ when $\frac{f(x)}{g(x)}\to 1$ as $x\to x_0$. 

Using l'H\^{o}pital's rule we see that $F_\nu(t)\approx K t^{-2s+1}$ as $t\to\-\infty$ and, as a consequence,  $F^{-1}_\nu(h)\approx K h^{1/(1-2s)}$ as $h\to 0+$. Furthermore,
\begin{equation}\label{Eq.Derivada1}
F_G'(h)= \frac{f_\nu(F_\nu^{-1}(h) + \mu)}{f_\nu(F_\nu^{-1}(h))}= \frac{(\nu +(F_\nu^{-1}(h) + \mu)^2)^{-s}}{(\nu + (F_\nu^{-1}(h))^2)^{-s}} \to 1, \quad \mbox{as } h\to 0+. 
\end{equation}

Some simple but tedious computations give that
\[
F_G''(h)\approx  K\frac{\mu \left(F_\nu^{-1}(h)\right)^{2s} +O\left(\left(F_\nu^{-1}(h)\right)^{2s-1}\right)}{\nu + (F_\nu^{-1}(h))^2};
\]
therefore, $F_G''(h)\approx K \left(F_\nu^{-1}(h)\right)^{2s-2}$. Consequently, $F_G''(h)\approx K h^{(2s-2)/(1-2s)}$. Now, applying 
l'H\^opital's rule twice we get that $\Delta(h)\approx K h^{\frac{2s-2}{1-2s}+2}=K h^{\frac{2s}{2s-1}}= K  h^{\frac{\nu +1}{\nu}}$, that is,
\[
\Delta(h) = K h^{(\nu+1)/\nu} + o(h^{(\nu+1)/\nu}) \mbox{ for some } K\neq 0 \mbox{ as } h\to 0+.
\]

We see from \eqref{Eq.Derivada1} that $F_G'$ is bounded. Hence $F_G$ is Lipschitz and Theorem \ref{Theo_0_1.1} can be applied here with $r_R=\frac{\nu+1}{\nu}$ to obtain 
$$(n+m)^{\frac{\nu}{\nu+2}}\ell_{n,m}^{0}\convw T_{\frac{\nu+1}{\nu},\frac{\nu+1}{\nu}}(0;K,K)$$ for any $\nu>0$.}
\FIN
\end{Ejem}

\begin{Ejem}\label{Ejem.Normal}
{\rm
Let $F$ (resp. $G$) be centered (resp. with mean $\mu>0$) normal distributions with common variance $\sigma^2$. Let $f$  denote the density function of $F$. Now, $F_G(t)=F( F^{-1}(t)+\mu)$, $t\in[0,1]$ and
\[
F_G'(t)=
\frac{f(F^{-1}(t)+\mu)}{f(F^{-1}(t))}= e^{-(2\mu  F^{-1}(t)+\mu^2)/2\sigma^2} \to 
\infty,  \mbox{ as } t \to 0+.
\]
{
This implies that $F_G$ is not Lipschitz in a neighbourhood of $0$. However, we can 
use the fact that $\ell_{n,m}^{t_0}=-\big(\int_0^\eta I(F_n^{-1}(t)\leq G_m^{-1}(t))dt-\int_0^\eta I(F^{-1}(t)\leq G^{-1}(t))dt\Big)=
-\big(\int_0^\eta I(F_n^{-1}(t)< G_m^{-1}(t))dt-\int_0^\eta I(F^{-1}(t)< G^{-1}(t))dt\Big)$. Hence, $\ell_{n,m}^{t_0}=-\tilde{\ell}_{m,n}^{t_0}$,
where $\tilde{\ell}_{m,n}^{t_0}$, as before, denotes the same statistic as $\ell_{n,m}^{t_0}$, but exchanging the roles of $F$ and $G$. Now $G_F'(0)=0$ and $0$ is an isolated {regular} contact point as in \eqref{Eq.DesarrolloDelta}, with $r_R=1$ and $C_R=-1$. Using Theorem \ref{Theo_0_1.1} we conclude that
$$(n+m)\ell_{n,m}^{t_0}\convw {\textstyle \lambda (1-\lambda) \int_0^\infty I_{\big\{- (1-  \lambda ) S^2_{\lceil \lambda  y \rceil} >0\big\} }dy=
\lambda (1-\lambda) \int_0^\infty I_{\big\{(1-  \lambda ) S^2_{\lceil \lambda  y \rceil} <0\big\} }dy}=0,$$
since, a.s.,  $S^2_i>0$ for every $i\geq 1$. Thus the rate of convergence in this example is faster than $(n+m)^{-1}$.
\FIN
}}
\end{Ejem}

We explore now some consequences of Theorem \ref{TCL.principal}. If the extremal contact points have a non-null contribution to the limiting distribution, then this cannot be normal. We pay now attention to obtaining conditions under which $\sqrt{n+m}\ell_{n,m}^0$ vanishes (and similarly for the upper extreme). The special attention to the rate $\sqrt{n+m}$ is due to the fact that it is the only one which can result in a normal limit. Of course, Theorem \ref{Theo_0_1.1} provides some answer to this problem, but we will give here simpler sufficient conditions.

If the supports of $F$ and $G$ are bounded and
\begin{equation}\label{tails}
\liminf  |F^{-1}(t)-G^{-1}(t)| > 0  \mbox{ when } t\to 0+ \mbox{ or $t\to 1-$} ,
\end{equation}
then $\ell_{n,m}^0$ and $\ell_{n,m}^1$
can be dealt with as in Lemma \ref{localizacion} to see that they eventually vanish. 

Note that, in the case of non-bounded support, (\ref{tails}) does not exclude  that 0 or 1 could be contact points (recall Example \ref{Ejem.Normal}). For this case the following  criterion on the tails can be useful to guarantee asymptotic negligibility of $\ell_{n,m}^0$ and $\ell_{n,m}^1$ in presence of inner contact points:
\begin{equation}\label{cond3}
\int_{(0,\varepsilon)\cup(1-\varepsilon,1)} \Big(\frac{\sqrt{t(1-t)}}{f(F^{-1}(t))} \Big)^p dt<\infty \mbox{ and } \ \int_{(0,\varepsilon)\cup(1-\varepsilon,1)}  \Big(\frac{\sqrt{t(1-t)}}{g(G^{-1}(t))} \Big)^p dt<\infty,
\end{equation}
for some $p>1$ and $\varepsilon>0$.

 {
In fact, let us assume that (\ref{cond3}) and (\ref{tails}) hold and that, for instance, $0$ is  a contact point and that $\inf (F^{-1}(t)-G^{-1}(t))>\delta>0$ on $(0,\eta)\subset (0,\varepsilon)$. We then focus   on the integral
$\int_{0}^\eta I_{\{F_n^{-1}(t)-G_m^{-1}(t)\leq 0\}} dt$, noting that 
$$(F^{-1}-G^{-1}>\delta)\cap(F_n^{-1}-G_m^{-1}\leq 0)\subset (|F_n^{-1}-F^{-1}|>\delta/2)\cup (|G_m^{-1}-G^{-1}|>\delta/2)$$

Now,  (\ref{cond3}) gives  that
\begin{eqnarray*}
\lefteqn{{n}^{1/2}\int_{0}^\eta I_{\{F_n^{-1}(t)-G_m^{-1}(t)\leq 0\}} dt} \\
&\leq &
{n}^{1/2}\left(\int_{0}^\eta I_{\{|F_n^{-1}(t)-F^{-1}(t)|\geq \delta/2\}} dt+\int_{0}^\eta I_{\{|G_m^{-1}(t)-G^{-1}(t)|\geq \delta/2\}} dt\right)\\
&\leq&  
\frac{n^{-\frac {p-1} {2}}}{(\delta/2)^p}\left(\int_{0}^\eta |\sqrt{n} (F_n^{-1}(t)-F^{-1}(t))|^{p} dt+\int_{0}^\eta |\sqrt{n} (G_m^{-1}(t)-G^{-1}(t))|^{p} dt\right)\convp 0,
\end{eqnarray*}
where the last convergence follows from the fact that by  (\ref{cond3}) and Theorem 5.3, p. 46 in \cite{BobkovLedoux}, 
the integrals in  parentheses are stochastically bounded. 
}
}

Now, we are ready for a general result for probabilities with smooth densities, $f$ and $g$. Assuming enough differentiability, we write $h(t)=F_G(t)-t$ (the function used to obtain \eqref{k.odd}) and, for any $k \in \Nat$,  define the sets
\[
\Gamma_k:= \left\{t \in \Gamma: h^{j)}(t)=0, j=0,\ldots,k-1 
\mbox{ and } h^{k)}(t) \neq 0\right\}.
\] 

Notice that the set $\Gamma_k$ is the set of contact points with intensity $k$ and let $k_0:=k_0^{F,G}= \sup \{k: \Gamma_k \neq \emptyset\}$. For points in $\Gamma_k$ the derivatives of $h$ can be easily related of the derivatives of $f$ and $g$, as follows.
\begin{Lemm}\label{LemaTontisimo}
If  $t_0 \in \Gamma_k$ for some $k\geq 1$,  and we denote $x_0=F^{-1}(t_0)$, then, 
\[
h^{k)}(t_0) 
=
\left\{
\begin{array}{ll}
\frac  {f(x_0)} {g(x_0)} -1, & \mbox{ if } k=1
\\
[2mm]
\frac {f^{k-1)}(x_0)-g^{k-1)}(x_0)} {f^k(x_0)}  & \mbox{ if } k>1,
\end{array}
\right.
\]
with $f(x_0) \neq g(x_0)$ in the first case and $f^{k-1)}(x_0)\neq g^{k-1)}(x_0)$ in the second one.
\end{Lemm}
\medskip

Combining \eqref{k.odd}, \eqref{1.odd}  and \eqref{k.even} with the above considerations we obtain the following version of Theorem \ref{TCL.principal}. 

\begin{Theo}\label{dosmuestras_2}
Assume that $F$ and $G$ have positive densities $f$ and $g$ on possibly unbounded intervals which are $k_0$ times continuously differentiable. Assume further  that the set of contact points is finite with maximal intensity $k_0$ and that  condition (\ref{tails})  holds. Suppose in addition that either the supports are bounded or that  condition (\ref{cond3}) is satisfied. 
Then, if $B_1$ and $B_2$  are independent Brownian bridges, and $n, m\to \infty$ with $\frac n {n+m} \to \lambda \in (0,1)$, 
\begin{enumerate}
\item[(i)] if $k_0=1$ and $x_i=F^{-1}(t_i)$,
$$(n+m)^{1/2}(\gamma(F_n,G_m)-\gamma(F,G)) \convw \sum_{t_i\in \Gamma_1}\Big({\textstyle  \frac{g(x_i)}{|f(x_i)-g(x_i)|}\frac {B_1(t_i)}{\sqrt{\lambda}}+\frac{f(x_i)}{|f(x_i)-g(x_i)|}\frac {B_2(t_i)}{\sqrt{1-\lambda}}}\Big),$$ 
\item[(ii)]
if $k_0\geq 3$ is odd
$$
(n+m)^{\frac 1 {2k_0}}(\gamma(F_n,G_m)-\gamma(F,G))\convw \sum_{t_i\in \Gamma_{k_0}}{\textstyle \Big(\frac{k_0!}{|h^{k_0)}(t_i)|}\Big)}^{1/k_0}\big(((B_\lambda(t_i))^{1/k_0})^+-((B_\lambda(t_0))^{1/k_0})^-\big),$$
\item[(iii)]
if $k_0$ is even
$$(n+m)^{\frac 1 {2k_0}}(\gamma(F_n,G_m)-\gamma(F,G))\convw \sum_{t_i\in \Gamma_{k_0}} \mbox{\em sgn}(h^{k_0)}(t_i))2{\textstyle \Big(\frac{k_0!}{|h^{k_0)}(t_i)|}}\Big)^{1/k_0}\big((B_\lambda(t_i))^{\mbox{\em \tiny sgn}( h^{k_0)}(t_i))}\big)^{1/k_0}. $$
\end{enumerate}
\end{Theo}

We see from Theorem \ref{dosmuestras_2} that asymptotic normality (arguably, the most useful case for statistical applications)
holds, with the standard $\sqrt{n+m}$ rate, only when $F$ and $G$ have a finite number of `simple' crossings.  In all the other cases we get a slower rate and a nonnormal limit.

While Theorem \ref{TCL.principal} (hence, also Theorem \ref{dosmuestras_2}) involves only the case when $\Gamma^*=\Gamma^*_F$ consists of regular contact points, the comments about virtual contact points between $F_G$ and the identity that led to \eqref{tangencia.virtual} apply to the global analysis of $\gamma(F_n,G_m)$. As an important example, 
we consider the case
when $F$ and $G$ are finitely supported. More precisely, let us assume $F$ and $G$ have a finite support $x_1<x_2<\dots<x_k$, with probabilities $p_1,p_2,\dots,p_k$ and $q_1,q_2,\dots,q_k$, respectively, with $p_i+q_i>0$ (although $p_i$ or $q_i$ could be null), $i=1,\ldots,k$. 
We set $P_i:=\sum_{j=1}^ip_j$ and $Q_i:=\sum_{j=1}^iq_j$, $i=1,\ldots k-1$. Then, $F_G(t)=P_i$ for $t\in (Q_{i-1},Q_i]$. Hence, the 
only possible inner contact points are $P_i, Q_i$, $i=1,\ldots,k-1$ and all the possible contact points are either horizontal crossings 
($P_i$ if $Q_{i-1}<P_i<Q_i$), vertical crossings ($Q_i$ if $Q_i<P_i<Q_{i+1}$), upper tangency points ($Q_i$ if $Q_{i-1}<Q_i=P_i<P_{i+1}$)
or lower tangency points ($P_i$ if $P_{i-1}<P_i=Q_{i-1}<Q_i$ ), using the same terms as in the discussion following Proposition \ref{Theo_Sing_disc}.
Combining that discussion with Corollary \ref{localizacion2} we obtain the following consequence.

\begin{Theo}\label{finite.support}
With the above notation, if  $F$ and $G$ are finitely supported and $\mathcal{H},\mathcal{V},\mathcal{U}$ and $\mathcal{L}$ denote, respectively, the sets of horizontal crossing, vertical crossing, upper tangency and lower tangency points for $F$ and $G$, then, assuming that $\frac{n}{n+m}\to \lambda \in(0,1)$,
\begin{eqnarray*}
\lefteqn{\sqrt{n+m}(\gamma(F_n,G_m)-\gamma(F,G))\convw \sum_{t\in \mathcal{H}}{\textstyle \frac{B_1(t)}{\sqrt{\lambda}}}-
\sum_{t\in \mathcal{V}}{\textstyle \frac{B_2(t)}{\sqrt{1-\lambda}}}}\hspace*{5cm}\\
&&+\sum_{t\in \mathcal{U}}{\textstyle \Big(\frac{B_1(t)}{\sqrt{\lambda}}-\frac{B_2(t)}{\sqrt{1-\lambda}}\Big)^+}-
\sum_{t\in \mathcal{L}}{\textstyle \Big(\frac{B_1(t)}{\sqrt{\lambda}}-\frac{B_2(t)}{\sqrt{1-\lambda}}\Big)^-},
\end{eqnarray*}
where $B_1$ and $B_2$ are independent Brownian bridges.
\end{Theo}

Similar to Theorem \ref{pieces}, we get a Gaussian limiting distribution only when all the contact points are crossing points (which, necessarily, have orders $r_L=r_R=1$). In the case $F=G$
we have $Q_{i-1}<Q_i=P_i<P_{i+1}$ for all $i$, that is, every $P_i$ is an upper tangency point and Theorem \ref{finite.support}
yields 
\begin{equation}\label{dos_ig2}
\sqrt{n+m} \gamma(F_n,G_m)\convw \sum_{i=1}^{k-1}{\textstyle \Big(\frac{B_1(P_i)}{\sqrt{\lambda}}-\frac{B_2(P_i)}{\sqrt{1-\lambda}}\Big)^+}.
\end{equation} 
Of course, using the fact that $\sqrt{1-\lambda}B_1-\sqrt{\lambda}B_2$ is a Brownian bridge, we can, equivalently, write \eqref{dos_ig2} as
$$\textstyle{\sqrt{\frac{nm}{n+m}}} \gamma(F_n,G_m)\convw \sum_{i=1}^{k-1} (B_{{1}}(P_i) )^+.$$

\bibliographystyle{abbrv}


\appendix

\section*{Appendix.}

\renewcommand{\thesubsection}{\Alph{subsection}}

\subsection{On the composite map $F_G$}\label{lips}
 
We collect here some useful fact about the transform $F_G$ (and $G_F$). At some points we have used the fact that, as a consequence of \eqref{compuesta}, for every measurable $A\subset[0,1]$,  
\begin{eqnarray} \nonumber
\ell \{t\in A  : t> F_G(t) \}&=&\ell \{t\in A  :  F^{-1}(t)> G^{-1}(t) \}, 
\\ \label{importante}
\ell \{t\in A : \mathbb U_n(t)>F_G(\mathbb V_m(t))\}&=& \ell \{t\in A  :  F^{-1}_n(t)> G^{-1}_m(t) \}.
\end{eqnarray} 
Looking at (\ref{importante}), the corresponding statement for $G_F$ would be 
\[
\ell \{t\in A : \mathbb V_m(t)>G_F(\mathbb U_n(t))\}= \ell \{t\in A  :  F^{-1}_n(t)< G^{-1}_m(t) \}.
\]
This shows that we can base our analysis indistinctly using $G_F$ or $F_G$, and, in particular, to study $\tilde{\ell}_{m,n}^{t_0}$ instead of $\ell_{n,m}^{t_0}$, (recall the discussion after Corollary \ref{localizacion2}) if 
\begin{itemize}
\item[i)]$\ell \{t\in (t_0-\eta,t_0+\eta)  :  F^{-1}(t)=G^{-1}(t) \}=0,$ 
\item[ii)]$P\left(\{ \ell \{t\in (t_0-\eta,t_0+\eta)  :  F^{-1}_n(t)=G^{-1}_m(t) \}>0\} \ \mbox{infinitely often} \right)=0,$
\item[iii)]$P\left( \ell \{t\in (t_0-\eta,t_0+\eta)  :  \mathbb V_m(t)=G_F(\mathbb U_n(t))\}>0\right)=0,$ and
\item[iv)]$P\left( \ell \{t\in (t_0-\eta,t_0+\eta)  :  \mathbb U_n(t)=F_G(\mathbb V_m(t))\}>0\right)=0$
\end{itemize}
hold. 
When $t_0$ is an isolated contact point then i) is satisfied. The other relations can be easily guaranteed taking into account the
next lemma and its  consequences.
\begin{Lemm}
{Let $X,Y$ be independent r.v.'s with respective d.f.'s $F$ and $G$.  Then $P(X=Y)=0$ if and only if $F$ and $G$ have no common discontinuity point.
}
\end{Lemm}

Therefore, if $F$ and $G$ have no common discontinuity point, the samples $\{X_1,\ldots,X_n\}$ and $\{Y_1,\ldots,Y_m\}$ are a.s. disjoint. Since these samples are the images of $F_n^{-1}$ and $G_m^{-1}$ respectively, the set $\{F_n^{-1}=G_m^{-1}\}$ must be a.s. empty. On the contrary, if there exists a common discontinuity point, $x$, for $F$ and $G$, then
\[
P\big(\ell\{F_n^{-1}=G_m^{-1}\}>0\big)
\geq
P\big(\ell\{F_n^{-1}=G_m^{-1}\}=1\big)
=P\big(X_1=x)^nP\big(Y_1=x)^m>0.
\]
This proves the following proposition.

\begin{Prop}\label{common.disc}
{Let $F_n^{-1}$ and $G_m^{-1}$ be the sample quantile functions based on independent samples of i.i.d. r.v.'s from the d.f.'s $F$ and $G$. Then $P(\ell \{F_n^{-1}=G_m^{-1})\}>0)>0$ for some $n,m$ if and only if $F$ and $G$ have a common discontinuity point.
}
\end{Prop}

Elaborating on the same ideas, it easily follows the following summarizing proposition.

\begin{Prop}\label{ultima}
Relations iii) and iv) above always hold. Moreover, if $F$ and $G$ do not have common discontinuity points on the set $[F^{-1}(t_0-\eta_0),F^{-1}(t_0+\eta_0)],$ then ii) holds for every  $\eta\in (0,\eta_0).$ 
\end{Prop}

Proposition A.25 in \cite{BobkovLedoux}  provides simple necessary and sufficient conditions under which a quantile function is Lipschitz. We exploit that characterization to give here necessary and sufficient conditions under which $F_G$ is Lipschitz.
\begin{Prop}\label{ultimabis}
The transform $F_G$ is Lipschitz if and only if 
\begin{equation}\label{lips1}
 F^{-1} \mbox{ is increasing  on $[F_G(0),F_G(1)]$ and  {\em supp}}(F)\cap (G^{-1}(0), G^{-1}(1))\subset \mbox{\em supp}(G)
 \end{equation}
and there exists some $\delta>0$ such that 
$$\lim \sup_{y\to x, y>x}\frac {G(F^{-1}(y)-)-G(F^{-1}(x)-)}{y-x}>\delta \mbox{ a.e. on {\em supp}}(F).$$ Condition \eqref{lips1} can be equivalently stated as
\begin{equation}\label{lips02} 
F \mbox{  is continuous  and } G \mbox{  increasing on  {\em supp}}(F)\cap (G^{-1}(0), G^{-1}(1)).
\end{equation}
\end{Prop}

\medskip
\noindent \textbf{Proof.}  We set $H^{-1}(t)=F_G(t-)$, $t\in(0,1)$ and note that $H^{-1}$ is left-continuous, hence a quantile function, and also that $H^{-1}$ is Lipschitz if and only if $F_G$ is Lipschitz. We write $H$ for the associated d.f., namely, $H(x)=\ell\{t:\, F_G(t-)\leq x\}=\ell\{t:\, F_G(t)\leq x\}$. By Proposition A.25 in \cite{BobkovLedoux} $H^{-1}$ is Lipschitz if and only if the associated probability is supported 
in a finite interval and its absolutely continuous component has a density separated from zero on that interval. 

Since the support of the law $\mathcal L(H^{-1})$ is contained in $[0,1]$, the first condition will hold if and only $H$ is strictly increasing on  $[F_G(0), F_G(1)]$ (see Proposition A.7 in \cite{BobkovLedoux}). This, in turn, is equivalent to \eqref{lips1} and to \eqref{lips02}.

In fact, to check the equivalence to (\ref{lips1}),  we note that if $H$ is increasing on $[F_G(0),F_G(1)]$ then for every $a,b \in (F_G(0),F_G(1))$, $a<b$, we have $\ell\{t: F_G(t)\in (a,b)\}>0$, which holds if and only if $\ell\{t:F_G(t)\in [a,b)\}>0$ for every $a,b \in (F_G(0),F_G(1))$, $a<b$, hence, if and only if $\ell\{t:G^{-1}(t)\in [F^{-1}(a),F^{-1}(b))\}>0$. This implies that $F^{-1}$ must be  increasing on  $(F_G(0),F_G(1)).$ 

Moreover, if $x\in (G^{-1}(0), G^{-1}(1))\setminus\mbox{supp}(G),$ then $x\in (G^{-1}(t^*),G^{-1}(t^*+))$ for some $t^*\in (0,1)$, while if $x\in \mbox{supp}(F)$, taking $\delta>0$ small enough such that $(x-\delta,x+\delta)\subset (G^{-1}(t^*),G^{-1}(t^*+))$, we would have $t^*=G(x-\delta)=G(x+\delta)$. Thus, $G^{-1}(t^*)<x-\delta<x+\delta<G^{-1}(t^*+)$, and $F(G^{-1}(t^*))\leq F(x-\delta)<F(x+\delta)\leq F(G^{-1}(t^*+))$. This would avoid the continuity (thus Lipschitzianity) of $F_G$. That (\ref{lips1}) implies that $H$ is increasing on $[F_G(0),F_G(1)]$ is immediate from the considerations above. The conditions in (\ref{lips02}) are easily seen to be equivalent to that in (\ref{lips1}).

For the condition involving the density of the absolutely continuous part, $F$ being continuous on $(G^{-1}(0), G^{-1}(1))$, $F_G$ is itself the quantile function of its distribution. The left-continuous version of H can be written as
$$
H(x-)=\ell\left\{t: F(G^{-1}(t))<x\right\}=\ell \left\{t: G^{-1}(t)<F^{-1}(x)\right\},
$$ 
and the density corresponding to the absolutely continuous part of $H$, should satisfy 
$$h(x)=\lim \sup_{y\to x, y>x}\frac {H(y)-H(x)}{y-x}>\delta \mbox{ a.e. on supp}(F)$$
for some $\delta>0$. This completes the proof. \FIN

\begin{Nota}\label{lips2} 
{\rm
Our analysis of the local behaviour of Galton's statistic around a contact point, $t_0$, required  $F_G$ to be Lipschitz on a neighbourhood $(t_0-\eta,t_0+\eta)$. This can be characterized in the same way as (\ref{lips1}) and (\ref{lips02}). Without trying to give the best possible result, this can be guaranteed, e.g., if   $G^{-1}$ is increasing and continuous on $(t_0-\eta,t_0+\eta)$ except perhaps at $t_0$, $F$ is continuous on $G^{-1}((t_0-\eta,t_0+\eta))$, and the derivatives of $F$ and $G$ (that exist almost everywhere) satisfy $\mbox{\em ess inf }\{\frac {G'(x)}{F'(x)}, x\in G^{-1}((t_0-\eta,t_0+\eta))\cap \mbox{supp}(F) \}>0$.
}
\end{Nota}

\subsection{Approximation of uniform quantile processes}

The following result has been used extensively in this paper. It is a consequence of a refined version of the Komlos-Major-Tusnady construction for the quantile process (see, e.g., Theorem 3.2.1, p. 152
in \cite{CsorgoHorvath}), from which we know that there exists a sequence of Brownian bridges on $[0,1]$, $\{B_{n}\}$, versions of $u_n$  and positive constants, $C_1,C_2$ and $C_3$, such that
\begin{equation} \label{approximation.KMT}
P\Big\{\sup_{0\leq t\leq 1} |u_n(t)-B_{n}(t)| >{\textstyle \frac{x+C_1 \log n}{\sqrt{n}}} \Big\}\leq C_2 e^{-C_3 x},\quad x>
0.
\end{equation}
Making use of this construction for both quantile processes and taking $x=\frac a {C_3}\log n$  with $a>1$ and $K=\frac a {C_3}+C_1>0$, we obtain useful independent sequences of Brownian bridges  $\{B_{n}^F\}$,  $\{B_{m}^G\}$  and versions of $u_n$ and $v_m$.
 
\begin{Theo}\label{aproxBridge}
With the previous notation, in a probability one set, the sequences $\{B_{n}^F\}$,  $\{B_{m}^G\}$,  $\{u_n\}$ and   $\{v_m\}$
eventually satisfy
\begin{equation}\label{eventualbound2}
\sup_{0\leq t\leq 1} |u_n(t)-B_{n}^F(t)|\leq {\textstyle K\frac{\log n}{\sqrt{n}}}
\quad 
\mbox{ and }
\quad
\sup_{0\leq t\leq 1} |v_m(t)-B_{m}^G(t)|\leq {\textstyle K\frac{\log m}{\sqrt{m}}}.
\end{equation}
\end{Theo}

\end{document}